\documentclass[12pt]{article}

\usepackage[margin=1.20in]{geometry}

\usepackage{amssymb,amsmath}
\usepackage{amsthm,enumitem}
\usepackage[nocompress]{cite}
\usepackage{hyperref}
\usepackage{mathrsfs}
\usepackage{textcomp}
\usepackage{xcolor}

\hypersetup{
    colorlinks,
    citecolor=black,
    filecolor=black,
    linkcolor=black,
    urlcolor=black
}

\usepackage{verbatim}

\newtheorem{Theorem}{Theorem}[section]
\newtheorem{Lemma}[Theorem]{Lemma}
\newtheorem{Proposition}[Theorem]{Proposition}

\newtheorem{Remark}[Theorem]{Remark}
\theoremstyle{definition}

\theoremstyle{plain}

\newcommand{\ud}{\mathrm{d}}
\newcommand{\R}{\mathbb{R}}
\newcommand{\Sn}{\mathbb{S}^{n-1}}

\newcommand{\abs}[1]{\lvert#1\rvert}

\begin{document}

\title{\textbf{Global classification and local asymptotics for a Monge--Amp\`ere obstacle problem}}
\author{
Tianling Jin,\footnote{T. Jin was partially supported by NSFC grant 12122120 and Hong Kong RGC grants GRF 16303822 and GRF 16303624.}\quad
Connor Mooney,\footnote{C. Mooney was partially supported by NSF grant DMS-2143668.}\quad
Aranya Sen,\footnote{A. Sen was partially supported by NSF grant DMS-2143668.}\\[4pt]
Xushan Tu,\quad
Jingang Xiong\footnote{J. Xiong was partially supported by NSFC grant 12325104.}}

\date{\today}

\maketitle

\begin{abstract}
Let $n\geq2$ and $0<q<n$. We prove that every nontrivial convex entire solution of
\[
\det D^2v=v^q,\quad v\geq0,\quad\text{in }\R^n
\]
with bounded (possibly empty) coincidence set $\{v=0\}$ is affine equivalent to a radial solution.  Together with the classification for unbounded coincidence sets by Jin--Tu--Xiong \cite{JTX2025}, this completes the classification of nontrivial nonnegative convex entire solutions. The proof uses an affine-invariant monotonicity formula extending the two-dimensional formula in Mooney--Sen \cite{MooneySen2026}. Combining this global classification with level-set parametrization and spectral analysis, we also obtain precise asymptotics of local solutions at isolated zeros, settling the remaining case, namely, a singleton free boundary, in the regularity theory for the Monge--Amp\`ere obstacle problem studied in \cite{JTX2025}.

\medskip



\end{abstract}


\section{Introduction}
For $n\geq2$ and $0<q<n$, we say that two convex functions $u_1$ and $u_2$ are $q$-affinely equivalent if
\begin{equation*}
u_1(x)=\abs{\det A}^{-2/(n-q)}u_2(Ax+b)
\end{equation*}
for some $A\in GL(n)$ and $b\in\R^n$.

For a convex function $v$ on a convex domain $\Omega$, its Monge--Amp\`ere measure is defined by $\mathscr M_v(E):=\abs{\partial v(E)}$ for Borel sets $E\subset\Omega$, where $\partial v(E)$ denotes the subgradient image of $E$. Throughout the paper, convex solutions to Monge--Amp\`ere equations are understood in the Alexandrov sense.

\begin{Theorem}\label{thm:complete-classification}
Let $n\geq2$ and $0<q<n$. Suppose that $v\not\equiv0$ is a convex solution of
\begin{equation}\label{eq:entire-solution}
\det D^2v=v^q,\qquad v\geq0,\qquad\text{in }\R^n.
\end{equation}
If $K:=\{v=0\}$ is bounded, then $K$ is either empty, a singleton, or an ellipsoid, and in each case, $v$ is $q$-affinely equivalent to a radial solution that is unique up to scaling.
\end{Theorem}
Together with the classification for unbounded coincidence sets by Jin--Tu--Xiong \cite{JTX2025}, this completes the classification of nontrivial convex entire solutions of \eqref{eq:entire-solution}. For $n=2$, Theorem~\ref{thm:complete-classification} was recently proved by Mooney--Sen \cite{MooneySen2026}. These radial profiles are described in \eqref{eq:radial-family}.

For \(q=0\), the corresponding obstacle problem is
\begin{equation}\label{eq:entire-solution00}
\det D^2v=\chi_{\{v>0\}},\qquad v\geq0,\qquad\text{in }\R^n.
\end{equation}
See Corollary~1.2 of Jin--Tu--Xiong \cite{JTX2025} for the classification of solutions to \eqref{eq:entire-solution00}.  In fact, when $|K|=0$, then the classification follows from the  J\"orgens--Calabi--Pogorelov theorem \cite{Jorgens1954,Calabi1958,Pogorelov1972}. When $0<|K|<\infty$, the classification follows from J\"orgens \cite{Jorgens1955} and Jin--Xiong \cite{JinXiong2016} through the Legendre transform. When $|K|=\infty$, the classification follows from Theorem 1.1 in Jin--Tu--Xiong \cite{JTX2025}. For $q\geq n$, nonzero convex solutions are positive by the strong maximum principle~in D\'iaz--D\'iaz \cite{Diaz2015}. Hence, the coincidence set $\{v=0\}$ is empty in this case.

The proof of Theorem~\ref{thm:complete-classification} relies on the following affine-invariant monotonicity formula, which was obtained by Mooney--Sen \cite{MooneySen2026} in two dimensions.

\begin{Theorem}\label{thm:main-monotonicity}
Let  $a<b$ be two real numbers, $f$ be a nonnegative bounded Borel measurable function, and $v$ be a convex function in a convex domain $\Omega$ satisfying
\[
\det D^2v=f(v)\quad\text{in } \{a<v<b\}.
\]
For $s\in \R$, set $K_s=\{v<s\}$. Assume that $v$ is strictly convex and $C^1$ in $\{a<v<b\}$, and that $0\in K_a \subset K_b\Subset\Omega$. Then
\begin{equation}\label{eq:Qn-introduction}
Q_n[v](s)=n\,\frac{\bigl(\partial_s\abs{K_s}\bigr)^n}{\abs{K_s}^{n-1}}\,\mathscr M_v(K_s)
\end{equation}
is affine-invariant, continuous, and nonincreasing on $(a,b)$. Moreover, if $Q_n[v]$ is constant on $(a,b)$, then $Q_n[v]\equiv n^{n+1}|B_1|^2$ and $\{K_s\}_{s\in (a,b)}$ form a concentric family of homothetic ellipsoids.
\end{Theorem}

Since the above monotonicity formula holds for general right-hand sides of the form $f(v)$, rigidity results such as Theorem~\ref{thm:complete-classification} can be extended to more general right-hand sides, see e.g. Remarks \ref{rem:round-level-set} and \ref{rem:rough-rhs}. Earlier, Collins--Tong \cite{CollinsTong2025} introduced a monotonicity formula for optimal transport between convex domains. However, the two monotonicity formulae are fundamentally different. The monotonicity formula in Theorem \ref{thm:main-monotonicity} is motivated by affine convex geometry and measures the failure of level sets to form a concentric family of homothetic ellipsoids. In contrast, the monotonicity formula in \cite{CollinsTong2025} is inspired by minimal/maximal surface theory and measures the failure of solutions to be homogeneous.
 
When using the monotonicity formula to prove Theorem \ref{thm:complete-classification}, we first obtain the upper bound $Q_n[v]\leq n^{n+1}|B_1|^2$ from the affine isoperimetric inequality. Affine normalization and compactness then produce an entire blow-down limit on which $Q_n$ is constant. The equality case of the monotonicity formula identifies this constant as $n^{n+1}|B_1|^2$. Monotonicity forces $Q_n[v]$ to have the same constant value, and the classification follows from the equality case.

The Monge--Amp\`ere obstacle problem
\begin{equation*}
\det D^2v=G(x,v(x),Dv(x))v^q, \quad v\ge 0,\quad\text{in }B_1,
\end{equation*}
emerges from the $L_p$-Minkowski problem \cite{Lutwak1993,LutwakYangZhang2004,CW2006}. Here $G$ is assumed to be bounded above and below by two positive constants. Away from the coincidence set $K=\{v=0\}$, regularity of the solutions follows from Caffarelli's localization and interior theory \cite{Caffarelli1990,Caffarelli1990W2p,Caffarelli1991}. If $v$ is allowed to vanish somewhere on $\partial B_1$, the free boundary $\partial K$ may have a non-strictly convex part, and $v$ can be merely Lipschitz at such points \cite[Remark~2.19]{JTX2025}. Jin--Tu--Xiong \cite{JTXSingular2025} further studied this part of the free boundary and established sharp dimension bounds for its flat portions. Strictly positive boundary data, however, rule out non-strictly convex points of the free boundary. For $|K|>0$, the regularity of the solution and of the strictly convex part of the free boundary has been established under suitable hypotheses through the works of Savin \cite{Savin2005}, G\'alvez--Jim\'enez--Mira \cite{GalvezJimenezMira2015}, Huang--Tang--Wang \cite{HuangTangWang2024}, Daskalopoulos--Lee \cite{DaskalopoulosLee2012}, Jin--Tu--Xiong \cite{JTX2025} and Han--Huang \cite{HanHuang2026}. The remaining regularity problem near the strictly convex part of the free boundary concerns the case $|K|=0$, in which $K$ must be a singleton if it has an exposed point in $B_1$ (see \cite[Lemma~2.15]{JTX2025}). Daskalopoulos--Lee \cite{DaskalopoulosLee2012} posed the regularity of the solution at isolated zeros as an open question, even in two dimensions. Below, we settle this remaining case by establishing precise asymptotics of the solution near an isolated zero.
 
Without loss of generality, we suppose $\{v=0\}=\{0\}$. By \cite[Theorem~1.7]{JTX2025}, $v\in C^{1,\alpha_0}$ near the origin for some $\alpha_0\in(0,1)$, and
\begin{equation*}
c|x|^{\frac{1+\alpha_0}{\alpha_0}}\le v(x)\le C|x|^{1+\alpha_0}\quad\mbox{near }0
\end{equation*}
for some positive constants $c$ and $C$. To obtain precise asymptotics near the isolated zero, we consider
\begin{equation}\label{eq:intro-structural-coefficient}
\det D^2v=g(x)v^q,\quad v\ge0,\quad\text{in }B_1,
\end{equation}
where $g$ is H\"older continuous near the origin and $g(0)=1$.
\begin{Theorem}\label{thm:intro-local-leading}
Let $0<q<n$ and $\lambda>0$. Suppose $\lambda\leq g\leq\lambda^{-1}$, $g$ is $C^\gamma$ near $0$ for some $\gamma\in(0,1)$ and $g(0)=1$. Let $v$ be a convex solution of \eqref{eq:intro-structural-coefficient} with $\{v=0\}=\{0\}$. Then there exists $A_\infty\in SL(n)$ such that
\begin{equation}\label{eq:intro-v-leading}
v(A_\infty^{-1}x)=c_*\abs{x}^{\frac{2n}{n-q}}+O(\abs{x}^{\frac{2n}{n-q}+\gamma})\qquad\text{as }x\to0,
\end{equation}
where $c_*=\frac{(n-q)^{(n+1)/(n-q)}}{(2n)^{n/(n-q)}(n+q)^{1/(n-q)}}.$
\end{Theorem}

A related model with an isolated degeneracy is
\[
\det D^2v(x)=|x|^\sigma,\quad v\ge0,\quad\text{in }B_1,\qquad \sigma>0.
\]
Daskalopoulos--Savin \cite{daskalopoulos2009monge} studied this equation in dimension two, where solutions may exhibit radial or nonradial behavior near the origin. In contrast, under the assumptions of Theorem~\ref{thm:intro-local-leading}, solutions of \eqref{eq:intro-structural-coefficient} are asymptotically radial after a fixed affine change of coordinates. To prove this, we first use affine normalization, the classification in Theorem~\ref{thm:complete-classification}, and the monotonicity formula in Theorem~\ref{thm:main-monotonicity} to show that every normalized blow-up limit at an isolated zero is affinely equivalent to the homogeneous radial solution. The remaining difficulty is that the normalizing affine transformations may depend on the scale, so this alone does not give an asymptotic expansion in fixed coordinates. Following the support-function approach of Choi--Choi--Kim \cite{CCKGrowth2025,CCKExistence2025,CCK2024}, we rewrite the equation as a first-order system and normalize the rescaled level sets in minimal mean width position \cite{GiannopoulosMilman2000}. Spectral analysis around the round stationary solution and the Merle--Zaag lemma \cite{MZ1998} then yield decay estimates that ensure convergence of the normalizing affine transformations.

The classification problems considered here and in the works of Jin--Tu--Xiong \cite{JTX2025} and Mooney--Sen \cite{MooneySen2026} have analogues for the classical obstacle problem
\begin{equation}\label{eq:cobstacle}
\Delta u=\chi_{\{u>0\}},\quad u\ge 0,\quad\text{in }\R^n,
\end{equation}
and for the Alt--Phillips problem \cite{AltPhillips1986},
\begin{equation}\label{eq:cobstacleAP}
\Delta u=\gamma u^{\gamma-1}\chi_{\{u>0\}},\quad u\ge 0,\quad\text{in }\R^n,\quad  1<\gamma<2.
\end{equation}
The classification of solutions to \eqref{eq:cobstacle} with bounded coincidence sets follows from the works of Dive \cite{Dive1931}, Lewy \cite{Lewy1979}, Sakai \cite{Sakai1981}, DiBenedetto--Friedman \cite{DiBenedettoFriedman1986}, and Friedman--Sakai \cite{FriedmanSakai1986}. The general case was resolved by Eberle--Shahgholian--Weiss \cite{EberleShahgholianWeiss2023} for $n\geq6$ and Eberle--Figalli--Weiss \cite{EberleFigalliWeiss2025} in all dimensions. By contrast, the classification of global solutions to the Alt--Phillips problem \eqref{eq:cobstacleAP} remains largely unresolved; see \cite{bonorino2001regularity,BonorinoEtAl2005,wu2022fully,dipierro2023classification,SavinYuCones2025,SavinYuConcentration2026} for related classification and regularity results, including other exponent ranges and fully nonlinear equations.

The paper is organized as follows. In Section~\ref{sec:monotonicity-derivation}, we prove the monotonicity formula. In Section~\ref{sec:global-classification}, we establish affine normalization and compactness, study blow-up and blow-down limits, and prove Theorem~\ref{thm:complete-classification}. In Section~\ref{sec:blow-up-normalization}, we analyze the rescaled support-function system and prove Theorem~\ref{thm:intro-local-leading}.

\section{Proof of the monotonicity formula}
\label{sec:monotonicity-derivation}

Let $K\subset\R^n$ be a smooth, uniformly convex body with $0\in\operatorname{int}K$. The Gauss map sends each $x\in\partial K$ to its outer unit normal $\nu_K(x)\in\Sn$, and its inverse $X_K$ sends $\theta\in\Sn$ to the unique boundary point with outer unit normal $\theta$. The support function is
\[
\mathbf h(\theta)=X_K(\theta)\cdot\theta=\sup_{x\in K}x\cdot\theta,
\]
which measures the distance from the origin to the tangent hyperplane of $\partial K$ at $X_K(\theta)$. The inverse Gauss parametrization is given by
\[
X_K(\theta)=\mathbf h(\theta)\theta+\nabla_{\Sn}\mathbf h(\theta)\in\partial K.
\]
Pulling back the boundary measure $\mathcal H^{n-1}|_{\partial K}$ via $X_K$ defines the surface-area measure $S_K$ on $\Sn$, with density
\[
\ud S_K(\theta)=\det\bigl(\nabla_{\Sn}^2\mathbf h(\theta)+\mathbf h(\theta)I\bigr)\,\ud\theta,
\]
where $\nabla_{\Sn}^2\mathbf h$ denotes the spherical Hessian with respect to the round metric on $\Sn$, $I$ is the identity on $T_\theta\Sn$, and $\ud\theta$ is the standard surface measure on $\Sn$. The determinant above is the Jacobian of $X_K$, which is the reciprocal of the Gauss curvature.

For a general convex body $K\subset\R^n$, let $\mathbf h_K(p)=\sup_{x\in K}x\cdot p$ for $p\in\R^n$ and write $\mathbf h=\mathbf h_K|_{\Sn}$. The outer unit normal $\nu_K$ is defined at $\mathcal H^{n-1}$-almost every point of $\partial K$. For any Borel set $A\subset\Sn$, $S_K(A)$ is the $\mathcal H^{n-1}$-measure of the set of boundary points whose outer unit normal lies in $A$. Writing the Lebesgue decomposition $S_K=S_K^a+S_K^s$ with respect to $\ud\theta$, where $S_K^a$ and $S_K^s$ are the absolutely continuous and singular parts of $S_K$, Alexandrov's theorem ensures that $\mathbf h$ is twice differentiable almost everywhere on $\Sn$, and the same determinant gives the density of $S_K^a$ almost everywhere \cite[(10.43)]{Schneider2014}. Finally, the divergence theorem yields
\begin{equation}\label{eq:volume-identities}
\abs K=\frac1n\int_{\Sn}\mathbf h\,\ud S_K.
\end{equation}

We use the following second-variation estimate for the volumes of sublevel sets of a convex function. For a smooth family of uniformly convex level sets, it follows from the local Brunn--Minkowski inequality \cite[Theorem~5.3]{Milman2025}.

\begin{Lemma}\label{lem:sublevel-volume-second-variation}
Let $v$ be a convex function on a convex domain $\Omega$, and suppose that $K_s=\{v<s\}$ is nonempty and compactly contained in $\Omega$ for every $s\in(a,b)$. Let $\mathbf h(\theta,s)$ be the support function of $\overline{K_s}$ and $\rho(\theta,s)$ the density of the absolutely continuous part of the surface-area measure $S_{K_s}$. Then $\abs{K_s}^{1/n}$ is concave in $s$, and
\[
\frac{\ud^2}{\ud s^2}\abs{K_s}^{1/n}\leq\frac{1}{n\abs{K_s}^{(n-1)/n}}\int_{\Sn}\mathbf h_{ss}\rho(\theta,s)\,\ud\theta\qquad\text{in }\mathcal D'(a,b),
\]
where $\mathbf h_{ss}$ is the almost-everywhere second derivative of the concave function $\mathbf h(\theta, \cdot)$ in $s$.
\end{Lemma}

\begin{proof}
Since $v$ is convex, we have
\[
(1-\lambda)K_{s_1}+\lambda K_{s_2}\subset K_{(1-\lambda)s_1+\lambda s_2}\qquad \mbox{for all }s_1,s_2\in(a,b),\ 0<\lambda<1.
\]
Thus $\mathbf h(\theta,\cdot)$ is concave, and the Brunn--Minkowski inequality implies that $s\mapsto\abs{K_s}^{1/n}$ is concave. For fixed $s,r\in(a,b)$, Minkowski's first inequality \cite[Theorem~7.2.1]{Schneider2014}, together with \eqref{eq:volume-identities}, gives
\[
n\abs{K_s}^{(n-1)/n}\bigl(\abs{K_r}^{1/n}-\abs{K_s}^{1/n}\bigr)\leq\int_{\Sn}\bigl(\mathbf h(\theta,r)-\mathbf h(\theta,s)\bigr)\,\ud S_{K_s}(\theta).
\]
Applying the preceding inequality with $r=s+\varepsilon$ and with $r=s-\varepsilon$, and then adding the two inequalities, we obtain
\[
n\abs{K_s}^{(n-1)/n}\delta_\varepsilon^2\bigl(\abs{K_s}^{1/n}\bigr)\leq\int_{\Sn}\delta_\varepsilon^2\mathbf h\,\ud S_{K_s}\leq\int_{\Sn}\delta_\varepsilon^2\mathbf h\,\rho(\theta,s)\,\ud\theta,
\]
where $\delta_\varepsilon^2 u(s) =\frac{u(s+\varepsilon)+u(s-\varepsilon)-2u(s)}{\varepsilon^2}$ denotes the symmetric second difference, and we used the concavity of $\mathbf h(\theta,\cdot)$ about $s$ in the second inequality. Now multiply this inequality by a nonnegative $\varphi\in C_c^\infty(a,b)$ and integrate in $s$. As $\varepsilon\to0$, the measures $\delta_\varepsilon^2(\abs{K_s}^{1/n})\,\ud s$ converge weakly on $(a,b)$ to the nonpositive measure $\frac{\ud^2}{\ud s^2}\abs{K_s}^{1/n}$, so we may pass to the limit on the left. On the right, $\delta_\varepsilon^2\mathbf h\leq0$ and $\delta_\varepsilon^2\mathbf h\to\mathbf h_{ss}$ for almost every $(\theta,s)$, so Fatou's lemma applied to the negatives of the integrands gives the asserted inequality.
\end{proof}

\begin{proof}[\textbf{Proof of Theorem~\ref{thm:main-monotonicity}}]
For $s\in(a,b)$, let $\mathbf h(\theta,s)=\mathbf h_{\overline{K_s}}(\theta)$ denote the support function of $\overline{K_s}$, and let $S_{K_s}$ denote the surface-area measure of $\overline{K_s}$. Since $\partial K_s$ is $C^1$ and strictly convex, each $\theta\in\Sn$ is the outer unit normal at a unique point $X(\theta,s)\in\partial K_s$. The point $X(\theta,s)$ depends continuously on $(\theta,s)$. By assumption $K_a\neq\emptyset$, the gradient of $v$ is nonzero on $\{a<v<b\}$, and
\[
\mathbf h_s(\theta,s)=\frac{1}{\abs{Dv(X(\theta,s))}} >0.
\]
Moreover, $\mathbf h_s$ is continuous in $(\theta,s)$.

By definition, $Q_n$ is affine invariant. Since $\mathbf h_s$ is jointly continuous, the first variation formula for volume \cite[Lemma~7.5.3]{Schneider2014} yields
\begin{equation*}
\frac{\ud}{\ud s}\abs{K_s}=\int_{\Sn}\mathbf h_s\,\ud S_{K_s}>0.
\end{equation*}
The joint continuity of $\mathbf h_s$ and the weak continuity of $S_{K_s}$ under Hausdorff convergence imply that $\abs{K_s}$ is $C^1$ in $s$. For $a<s_1<s_2<b$, the Monge--Amp\`ere equation then implies
\[
\mathscr M_v(K_{s_2})-\mathscr M_v(K_{s_1})=\int_{\{s_1<v<s_2\}}f(v)\,\ud x=\int_{s_1}^{s_2}f(s)\frac{\ud}{\ud s}\abs{K_s}\,\ud s.
\]
Consequently, $\mathscr M_v(K_s)$ is locally absolutely continuous in $s$, and $Q_n[v]$ is continuous on $(a,b)$. Moreover, since $v$ is strictly convex, $\mathscr M_v(K_s)$ is strictly increasing in $s$. The coarea formula also shows that $\{f>0\}$ has positive measure in every open subinterval of $(a,b)$.

Write $w=\mathbf h_s$. The subgradient image $\partial v(K_s)$ is a star-shaped domain with radial function $w(\cdot,s)^{-1}$, so integration in polar coordinates gives
\begin{equation*}
n\mathscr M_v(K_s)=n\abs{\partial v(K_s)}=\int_{\Sn}w(\theta,s)^{-n}\,\ud\theta.
\end{equation*}
Moreover, for $a<s_1<s_2<b$ and a Borel set $A\subset\Sn$, let $E_A$ be the set of points $x$ satisfying $s_1<v(x)<s_2$ and $Dv(x)/\abs{Dv(x)}\in A$.  Computing $\mathscr M_v(E_A)$ from its subgradient image and then using the Monge--Amp\`ere equation gives
\begin{equation}\label{eq:mass-EA}
\begin{aligned}
\frac1n\int_A\bigl(w(\theta,s_2)^{-n}-w(\theta,s_1)^{-n}\bigr)\,\ud\theta&=\mathscr M_v(E_A)=\int_{s_1}^{s_2}f(s)\int_A w\,\ud S_{K_s}\,\ud s.
\end{aligned}
\end{equation}

Let $\rho(\theta,s)$ be the density of the absolutely continuous part of the surface-area measure $S_{K_s}$ with respect to $\ud\theta$. First, we assume that $v$ is smooth and uniformly convex and that $f$ is positive. Then $\ud S_{K_s}=\rho\,\ud\theta$, and \eqref{eq:mass-EA} yields
\[
w_s=-f(s)\rho(\theta,s)w^{n+2}.
\]
Differentiating $Q_n[v]$ and using Lemma~\ref{lem:sublevel-volume-second-variation} for the second derivative of $\abs{K_s}^{1/n}$, we obtain
\begin{align}
\frac{\ud}{\ud s}\log Q_n[v]&=\frac{f(s)}{\mathscr M_v(K_s)}\frac{\ud}{\ud s}\abs{K_s}+\frac{n}{\dfrac{\ud}{\ud s}\abs{K_s}^{1/n}}\frac{\ud^2}{\ud s^2}\abs{K_s}^{1/n}\label{eq:Qn-log-derivative}\\
&\leq\frac{f(s)}{\mathscr M_v(K_s)}\frac{\ud}{\ud s}\abs{K_s}-\frac{nf(s)}{\dfrac{\ud}{\ud s}\abs{K_s}}\int_{\Sn}\rho^2w^{n+2}\,\ud\theta\label{eq:Qn-log-derivative1}\\
&\leq 0,\label{eq:Qn-log-derivative2}
\end{align}
where the last inequality uses
\begin{equation}\label{eq:Qn-cauchy-schwarz}
\left(\frac{\ud}{\ud s}\abs{K_s}\right)^2=\left(\int_{\Sn}w\rho\,\ud\theta\right)^2\leq n\mathscr M_v(K_s)\int_{\Sn}\rho^2w^{n+2}\,\ud\theta.
\end{equation}

In the general case, we first show that the same equation for $w_s$ holds almost everywhere under the stated hypotheses. For almost every $s$ with $f(s)>0$, the Monge--Amp\`ere equation and the coarea formula imply that $D^2v>0$ at $\mathcal H^{n-1}$-almost every point of $\partial K_s$. At such points, the second fundamental form is the restriction of $D^2v/\abs{Dv}$ to the tangent space, hence the generalized Gauss curvature is positive. It follows from Hug~\cite[Theorem~3.7]{Hug1999} that
\[
\ud S_{K_s}=\rho(\theta,s)\,\ud\theta\qquad\text{for almost every }s\text{ with }f(s)>0.
\]
At these heights, the first variation formula reads $ \frac{\ud}{\ud s}\abs{K_s} =\int_{\Sn}w\rho\,\ud\theta$, and \eqref{eq:Qn-cauchy-schwarz} follows.
 
Moreover, substituting $\ud S_{K_s}=\rho(\theta,s)\,\ud\theta$ into \eqref{eq:mass-EA} and using Fubini's theorem and the arbitrariness of $A$, we find that $w(\theta,\cdot)$ is locally absolutely continuous for almost every $\theta$ and satisfies
\begin{equation}\label{eq:weak-support-equation}
w_s=-f(s)\rho(\theta,s)w^{n+2}\qquad\text{a.e. on }\Sn\times(a,b).
\end{equation}
Applying Lemma~\ref{lem:sublevel-volume-second-variation} and substituting $-\mathbf h_{ss}=-w_s=f(s)\rho(\theta,s)w^{n+2}$ from \eqref{eq:weak-support-equation}, we obtain
\begin{equation}\label{eq:Qn-volume-derivatives}
-n\abs{K_s}^{(n-1)/n}\frac{\ud^2}{\ud s^2}\abs{K_s}^{1/n}\geq\int_{\Sn}f(s)\rho(\theta,s)^2w^{n+2}\,\ud\theta\qquad\text{in }\mathcal D'(a,b).
\end{equation}
Note that \eqref{eq:Qn-volume-derivatives} holds automatically on $\{f=0\}$. Since the derivative of $\abs{K_s}^{1/n}$ is positive, continuous, and locally of bounded variation, the chain rule for continuous functions of bounded variation and the coarea formula give \eqref{eq:Qn-log-derivative} in $\mathcal D'(a,b)$, and the singular part of the distributional derivative of $\log Q_n[v]$ is nonpositive. Multiplying \eqref{eq:Qn-volume-derivatives} by $n/\bigl(\frac{\ud}{\ud s}\abs{K_s}\bigr)$ proves \eqref{eq:Qn-log-derivative1}. By \eqref{eq:Qn-cauchy-schwarz} a.e. on $\{f>0\}$, we obtain \eqref{eq:Qn-log-derivative2} a.e. on $\{f>0\}$. On $\{f=0\}$, both terms on the right-hand side of \eqref{eq:Qn-log-derivative1} vanishes. Thus $\frac{\ud}{\ud s}\log Q_n[v]\leq0$ in $\mathcal D'(a,b)$, and $Q_n[v]$ is nonincreasing.

Finally, suppose that $Q_n$ is constant on an interval. The preceding bounds force equality in \eqref{eq:Qn-cauchy-schwarz} for almost every $s$ in this interval with $f(s)>0$, so
\[
w^{n+1}\rho(\theta,s)=\frac{1}{n\mathscr M_v(K_s)}\frac{\ud}{\ud s}\abs{K_s}\qquad\text{a.e. where }f(s)>0.
\]
On $\{f=0\}$, we have $\bigl(\mathscr M_v(K_s)\bigr)'=0$ and $w_s=0$. Thus \eqref{eq:weak-support-equation} and $\bigl(\mathscr M_v(K_s)\bigr)'=f(s)\frac{\ud}{\ud s}\abs{K_s}$ imply
\[
\partial_s\bigl(\mathscr M_v(K_s)^{1/n}w\bigr)=0\qquad\text{a.e. in the interval}.
\]
Since $Q_n=n^{n+1}\mathscr M_v(K_s)\bigl(\frac{\ud}{\ud s}\abs{K_s}^{1/n}\bigr)^n$, for each $\theta$ the support function $\mathbf h(\theta,s)$ is affine as a function of $\abs{K_s}^{1/n}$. The equality case of the Brunn--Minkowski inequality therefore gives $K_s=z+R(s)E$, where $z$ is fixed and $E$ is a convex body. By the first variation formula, $R'(s)>0$, so $w=R'(s)\mathbf h_E>0$ implies that $0\in\operatorname{int}E$. Choose a height $s$ with $f(s)>0$ such that equality holds in \eqref{eq:Qn-cauchy-schwarz} and $\ud S_{K_s}=\rho\,\ud\theta$. Substituting $K_s=z+R(s)E$ into the identity for $w^{n+1}\rho$ gives
\[
\mathbf h_E^{n+1}\,\ud S_E=C\,\ud\theta, \qquad \text{for some }C>0.
\]
By the theorem of Blaschke and Deicke, closed convex hypersurfaces of constant centro-affine curvature are origin-centered ellipsoids; see also \cite[Proposition~1.2 and p.~36]{CW2006}. Therefore, $E$ is an ellipsoid, and affine invariance reduces the calculation to $K_s=B_{R(s)}$, where $\omega_n=|B_1|$, $\abs{K_s}=\omega_nR^n$, and $\mathscr M_v(K_s)=\omega_n(R')^{-n}$. In this case, $Q_n=c_n=n^{n+1}\omega_n^2$.
\end{proof}

We also recall the affine surface area
\[
\operatorname{as}(K)=\int_{\Sn}\det\bigl(\nabla_{\Sn}^2\mathbf{h}+\mathbf{h}I\bigr)^{n/(n+1)}\,\ud\theta;
\]
see \cite[(10.34)]{Schneider2014}. The affine isoperimetric inequality of Blaschke \cite[Theorem~4.6]{Lutwak1991} states that
\begin{equation}\label{eq:affine-isoperimetric-inequality}
\operatorname{as}(K)^{n+1}\leq c_n\abs K^{n-1},\qquad c_n=n^{n+1}\omega_n^2,\qquad \omega_n=\abs{B_1},
\end{equation}
with equality if and only if $K$ is an ellipsoid.

\begin{Lemma}\label{lem:Qn-upper-bound}
Let $f$ be a nonnegative Borel measurable function, and let $v\in C^1(\Omega)$ be a nonnegative convex solution of
\[
\det D^2v=f(v).
\]
Suppose that $0\in K_{s_0}\Subset\Omega$ for some $s_0>\min_\Omega v$ and that $f\in L^\infty_{\mathrm{loc}}((\min_\Omega v,s_0))$. If $v$ is strictly convex in $\Omega\setminus\{v=\min_\Omega v\}$, then
\[
Q_n[v](s)\leq c_n\qquad\text{on }(\min_\Omega v,s_0).
\]
\end{Lemma}
\begin{proof} 
By translation, we may assume that $v(0)=\min_\Omega v$. By the definition of $Q_n[v]$ in \eqref{eq:Qn-introduction},
\[
\frac{\mathscr M_v(K_s)}{Q_n[v](s)}=\frac{1}{n^{n+1}\bigl(\frac{\ud}{\ud s}\abs{K_s}^{1/n}\bigr)^n}.
\]
On compact subintervals of $(\min_\Omega v,s_0)$, the chain rule and \eqref{eq:Qn-volume-derivatives} give the distributional inequality
\[
\begin{aligned}
\frac{\ud}{\ud s}\left(\frac{\mathscr M_v(K_s)}{Q_n[v](s)}\right)&=-\frac{\dfrac{\ud^2}{\ud s^2}\abs{K_s}^{1/n}}{n^n\bigl(\frac{\ud}{\ud s}\abs{K_s}^{1/n}\bigr)^{n+1}}\\
&\geq\frac{\abs{K_s}^{n-1}}{\bigl(\frac{\ud}{\ud s}\abs{K_s}\bigr)^{n+1}}\int_{\Sn}f(s)\rho^2w^{n+2}\,\ud\theta\\
&\ge 0.
\end{aligned}
\]
Since $\frac{\ud}{\ud s}\abs{K_s}=\int_{\Sn}w\rho\,\ud\theta$ for almost every $s$ with $f(s)>0$, the affine isoperimetric inequality \eqref{eq:affine-isoperimetric-inequality} yields
\[
\left(\frac{\ud}{\ud s}\abs{K_s}\right)^{n+2}\leq\operatorname{as}(K_s)^{n+1}\int_{\Sn}\rho^2w^{n+2}\,\ud\theta\leq c_n\abs{K_s}^{n-1}\int_{\Sn}\rho^2w^{n+2}\,\ud\theta.
\]
Thus,
\[
\frac{\ud}{\ud s}\left(\frac{\mathscr M_v(K_s)}{Q_n[v](s)}\right)\geq\frac{f(s)}{c_n}\frac{\ud}{\ud s}\abs{K_s}=\frac{1}{c_n}\frac{\ud}{\ud s}\mathscr M_v(K_s).
\]
For $\min_\Omega v<a<s<s_0$, integrating from $a$ to $s$ gives
\[
\frac{\mathscr M_v(K_s)}{Q_n(s)}\geq\frac{\mathscr M_v(K_a)}{Q_n(a)}+\frac{\mathscr M_v(K_s)-\mathscr M_v(K_a)}{c_n}\geq\frac{\mathscr M_v(K_s)-\mathscr M_v(K_a)}{c_n}.
\]
Since $v\in C^1(\Omega)$, we have $\mathscr M_v(\{v=\min_\Omega v\})=0$, and hence $\mathscr M_v(K_a)\to0$ as $a\to \min_\Omega v$. Passing to the limit in the preceding inequality yields $Q_n[v](s)\le c_n$.
\end{proof}

\begin{Remark}\label{lem:positive-level-convergence-to-core}
\label{lem:Qn-free-boundary-limit}
Let $0<q<n$, and let $v$ be a nonnegative convex solution of
\begin{equation}\label{eq:global-obstacle}
\det D^2v=v^q, \quad v\ge 0,\quad\text{in }\Omega.
\end{equation}
Suppose that $K_{s_0}\Subset\Omega$ for some $s_0>0$ and that $\abs K>0$. By \cite[Theorem~6.1]{JTX2025} and Legendre duality, $\partial K$ is $C^2$ and uniformly convex, and $\mathbf h_{K_s}\to\mathbf h_K$ in $C^2(\Sn)$ as $s\to0$. Integrating \eqref{eq:weak-support-equation} and using this convergence gives
\begin{equation*}
\lim_{s\to0}Q_n[v](s)=\frac{\operatorname{as}(K)^{n+1}}{\abs K^{n-1}}.
\end{equation*}
\end{Remark}

\begin{Remark}\label{rem:round-level-set}
We obtain the following generalization of Mooney--Sen \cite[Theorem~6.1]{MooneySen2026}. Let $n\geq2$, and let $v\in C^1(B_R)\cap C(\overline{B_R})$ be a nonnegative convex solution of
\[
\det D^2v=f(v)\quad\text{in }B_R,\qquad v=1\quad\text{on }\partial B_R,
\]
where $f$ is Borel measurable, with $f\in L^\infty_{\mathrm{loc}}((\min_{B_R}v,1])$. Suppose that $v$ is strictly convex in $\{v>\min_{B_R}v\}$. Then every sublevel set $K_s$, $\min_{B_R}v<s<1$, is a ball centered at $0$.

Indeed, Lemma~\ref{lem:Qn-upper-bound} gives $Q_n[v](s)\leq c_n$. Since $K_s\to B_R$ as $s\to1^-$, concavity of $\abs{K_s}^{1/n}$ and Minkowski's first inequality give
\[
n\abs{B_1}^{(n-1)/n}\frac{\ud}{\ud s}\abs{K_s}^{1/n}\geq\frac{1}{1-s}\int_{\Sn}\bigl(R-\mathbf h(\theta,s)\bigr)\,\ud\theta.
\]
The right-hand side is the average of the nonincreasing integral $\int_{\Sn}w(\theta,r)\,\ud\theta$ over $s<r<1$, so it has the same limit as $\int_{\Sn}w(\theta,s)\,\ud\theta$ as $s\to1^-$. Since $f$ is bounded near $1$, $n\mathscr M_v(K_s)=\int_{\Sn}w^{-n}\,\ud\theta$ has a finite limit. Consequently,
\[
\lim_{s\to1^-}Q_n[v](s)\geq\frac{1}{\abs{B_1}^{n-1}}\lim_{s\to1^-}\left(\int_{\Sn}w^{-n}\,\ud\theta\right)\left(\int_{\Sn}w\,\ud\theta\right)^n\geq c_n.
\]
Monotonicity and the equality case in Theorem~\ref{thm:main-monotonicity} yield the conclusion.
\end{Remark}

\section{Global classification}
\label{sec:global-classification}

\subsection{Preliminaries on affine normalization and compactness}

We establish the affine normalization and compactness needed for the blow-up and blow-down arguments.

Let $0<q<n$, let $\Omega$ be a convex domain, and let $v$ be a convex solution of
\begin{equation}\label{eq:local-composite-equation}
\det D^2v=g(x)v^q,\quad v\ge 0,\quad 0<\lambda\leq g\leq\lambda^{-1}\quad\text{in }\Omega.
\end{equation}
For $A\in SL(n)$ and $b\in\R^n$, set $T_{A,b}(x)=A(x-b)$. Set
\[
\alpha=\frac{n-q}{2n}.
\]
For $s>0$ and $T=T_{A,b}$, define
\[
v_{s,T}(y)=s^{-1}v\bigl(T^{-1}(s^\alpha y)\bigr)=s^{-1}v\bigl(b+s^\alpha A^{-1}y\bigr),\qquad y\in\Omega_{s,T}:=s^{-\alpha}T(\Omega).
\]
Then
\[
K_r(v_{s,T})=s^{-\alpha}T(K_{rs}(v))=s^{-\alpha}A\bigl(K_{rs}(v)-b\bigr)\qquad (r>0),
\]
and equation \eqref{eq:global-obstacle} is invariant under this rescaling. The monotonicity functional satisfies
\begin{equation}\label{eq:Qn-scaling}
Q_n[v_{s,T}](r)=Q_n[v](sr).
\end{equation}
If $v$ solves \eqref{eq:local-composite-equation}, then $v_{s,T}$ satisfies
\[
\det D^2v_{s,T}=g_{s,T}(y)v_{s,T}^q\quad\text{in }\Omega_{s,T},\qquad g_{s,T}(y)=g\bigl(b+s^\alpha A^{-1}y\bigr).
\]
We write $v_{s,A}:=v_{s,T_{A,0}}$. If $v(b)=\min_\Omega v$, then $ v_{s,T}(0)=s^{-1}v(b)=\min_{\Omega_{s,T}}v_{s,T}$.

Throughout this subsection, we assume that
\[
v(0)=\min_\Omega v<1,\qquad 0\in K_1\Subset\Omega,
\]
and set $K=\{v=0\}$ and
\[
a_0(v)=\abs K^{\frac{2}{n-q}}+\min_\Omega v.
\]

\begin{Lemma}
\label{lem:section-normalization}
There exist $\varepsilon_0,c,C>0$, depending only on $n,q,\lambda$, such that, if $a_0(v)\leq\varepsilon_0$, there is a transformation $T=T_{A,b}$ with $v(b)=\min_\Omega v$ and $A\in SL(n)$ for which
\[
B_c\subset T(K_1)\subset B_C.
\]
\end{Lemma}

\begin{proof}
Fix $b\in \{v= \min v\}$ and translate coordinates so that $b=0$. Set $\delta=1-\min v$. Decreasing $\varepsilon_0$ so that $\varepsilon_0\leq1/2$, we have $\delta\geq1/2$, because $\min v=0$ when $K\neq\varnothing$ and $\min v=a_0(v)$ when $K=\varnothing$. The set $K_1$ is the section of $v$ at the origin with height $\delta$, and $g\,v^q\chi_{\{v>0\}}\leq\lambda^{-1}$ on $K_1$. The section estimate in \cite[Lemma~2.11]{JTX2025} gives
\[
\abs{K_1}\,\mathscr M_v(K_1)\geq c\delta^n.
\]
Since $\mathscr M_v(K_1)\leq\lambda^{-1}\abs{K_1}$, this gives $\abs{K_1}\geq c$.

By John's lemma, there are $r>0$ and an affine map $\widetilde T(x)=\widetilde A x+x_0$, with $\det\widetilde A=1$, such that
\[
B_r\subset \widetilde T(K_1)\subset B_{nr},\qquad r\asymp_n\abs{K_1}^{\frac{1}{n}}.
\]
Replacing $v$ by $v\circ\widetilde T^{-1}$ and retaining the same notation, we may assume that
\[
B_r\subset K_1\subset B_{nr},\qquad r\geq c.
\]
Fix $c=c(n,q)>0$ sufficiently small. If $K\neq\varnothing$ and $K$ intersected every ball $B_{cr}(y)$ with
\[
y\in\left\{\pm\frac r2 e_i:1\leq i\leq n\right\},
\]
then convexity would give $\abs K\geq cr^n$, contradicting $\abs K=a_0(v)^{(n-q)/2}$ after a further reduction of $\varepsilon_0$. Thus, whether or not $K$ is empty, one of these balls is contained in $K_1\setminus K$. On this ball, $v>0$, and $v\leq1$ on its boundary. By \cite[Lemma~2.5]{JTX2025}, we obtain
\[
cr^n\leq\abs{B_{cr}(y)}\leq C.
\]
This gives $r\leq C(n,q,\lambda)$.

The chosen minimum point of the transformed function is $x_0=\widetilde T(0)$. To remove the translation $x_0$, it remains to show that $\operatorname{dist}(x_0,\partial K_1)\geq c$. Since $\det\widetilde A=1$ and $v\leq1$ on $K_1$, we have $\mathscr M_v(K_1)\leq\lambda^{-1}\abs{K_1}\leq C$. Moreover, $1-v(x_0)\geq1/2$ for $a_0(v)$ sufficiently small. The Alexandrov maximum principle applied to $v-1$ at $x_0$ gives
\[
(1-v(x_0))^n\leq C_n\operatorname{dist}(x_0,\partial K_1)\operatorname{diam}(K_1)^{n-1}\mathscr M_v(K_1).
\]
Hence $\operatorname{dist}(x_0,\partial K_1)\geq c$. Since $K_1\subset B_{nr}$ and $r\leq C$, it follows that
\[
B_c\subset K_1-x_0\subset B_C.
\]
In terms of the original sublevel set, the translated set is $\widetilde A(K_1-b)$. Thus $T_{\widetilde A,b}(x)=\widetilde A(x-b)$ sends the chosen minimum point $b$ to the origin and has the required property.
\end{proof}

We say that $v$ is normalized if
\[
v(0)=\min_\Omega v<1,\qquad B_c\subset K_1(v)\subset B_C,\qquad K_1(v)\Subset\Omega,
\]
where $c,C$ are the fixed constants in Lemma~\ref{lem:section-normalization}. At height $s>0$, we call $T=T_{A,b}$ a normalizing affine transformation if $v(b)=\min_\Omega v$ and $v_{s,T}$ is normalized on $\Omega_{s,T}$. We then call $v_{s,T}$ a normalized rescaling and $A$ a normalizing matrix.

For $T=T_{A,b}$ with $v(b)=\min_\Omega v$, the scaling identities give
\[
a_0(v_{s,T})=s^{-1}a_0(v).
\]
Thus, if $\min_\Omega v<s$, $K_s(v)\Subset\Omega$, and $s\geq\varepsilon_0^{-1}a_0(v)$, Lemma~\ref{lem:section-normalization} applied to $v_{s,T_{I,b}}$ yields a normalizing affine transformation at height $s$.

\begin{Lemma}
\label{lem:coincidence-set-regularity}
\label{lem:positive-level-geometry}
There exists $\alpha_0=\alpha_0(n,q,\lambda)\in(0,1)$ such that $v\in C^{1,\alpha_0}(K_1)$ and $Dv=0$ on $K$. Suppose in addition that $g\in C_{\rm loc}^{k,\gamma}(\{\min_\Omega v<v<1\})$ for some integer $k\geq0$ and $\gamma\in(0,1)$. Then $v$ is locally uniformly convex and belongs to $C_{\rm loc}^{k+2,\gamma}$ in $\{\min_\Omega v<v<1\}$.

More precisely, suppose that, for some $\widetilde c,\widetilde C>0$,
\[
B_{\widetilde c}(0)\subset K_1\subset B_{\widetilde C}(0).
\]
Then, for every $v(0)<s_-<s_+<1$, there is a constant $C_*>0$, depending only on $n,q,\lambda,v(0),\widetilde c,\widetilde C,s_-,s_+$, such that
\[
\|v\|_{C^{1,\alpha_0}(\{s_-\leq v\leq s_+\})}\leq C_*.
\]
Furthermore, there is a constant $\lambda_*>0$, depending only on the same quantities and $\|g\|_{C^{k,\gamma}}$, such that
\[
\lambda_*I\leq D^2v\leq \lambda_*^{-1}I\quad\text{on }\{s_-\leq v\leq s_+\},\qquad\|v\|_{C^{k+2,\gamma}(\{s_-\leq v\leq s_+\})}\leq\lambda_*^{-1}.
\]
\end{Lemma}

\begin{proof}
The first assertion follows from the arguments in \cite[Proposition~2.18 and Theorems~1.4, 1.7]{JTX2025}. The Alexandrov estimate and convexity imply that $v$ is uniformly Lipschitz on $\{s_-\leq v\leq s_+\}$. We next establish the remaining estimates.

Let $\widehat s_-=(v(0)+s_-)/2$ and $\widehat s_+=(s_++1)/2$. For $x\in\{s_-\leq v\leq s_+\}$ and $h>0$, define
\[
S_h^v(x)=\{y\in\Omega:v(y)<v(x)+Dv(x)\cdot(y-x)+h\}.
\]
We claim that there is an $h_*>0$, depending only on $n,q,\lambda,v(0),s_-,s_+,\widetilde c$, and $\widetilde C$, such that
\[
S_{h_*}^v(x)\Subset K_{\widehat s_+}\setminus\overline{K_{\widehat s_-}}\qquad\text{for every }x\in\{s_-\leq v\leq s_+\}.
\]
Otherwise, there would be solutions $v_j$ satisfying $B_{\widetilde c}(0)\subset K_1(v_j)\subset B_{\widetilde C}(0)$ and points $x_j\in\{s_-\leq v_j\leq s_+\}$ such that the preceding compact containment fails at heights tending to zero. After passing to a subsequence, the Alexandrov estimate gives local uniform convergence of $v_j$ to a solution $v_\infty$; moreover, $x_j\to x_\infty\in\{s_-\leq v_\infty\leq s_+\}$. The vanishing section heights yield a nontrivial segment contained in $\{\widehat s_-\leq v_\infty\leq\widehat s_+\}$ on which $v_\infty$ is affine. This contradicts the strict convexity of $v_\infty$ in $\{\min v_\infty<v_\infty<1\}$, which follows from \cite[Propositions~2.18 and~3.1]{JTX2025} and the constant boundary value.

On $K_{\widehat s_+}\setminus\overline{K_{\widehat s_-}}$, the density $g v^q$ is bounded above and below by positive constants and $v$ is Lipschitz. After decreasing $h_*$ if necessary, the standard section estimates give uniform two-sided size bounds for $S_{h_*}^v(x)$. The required $C^{1,\alpha_0}$ and $C^{k+2,\gamma}$ estimates follow from the interior theory of Caffarelli \cite{Caffarelli1990,Caffarelli1991}, followed by Schauder estimates in $S_{h_*}^v(x)$.
\end{proof}

We next compare sublevel sets at different heights.

\begin{Lemma}
\label{lem:small-a0-annular-separation}
After decreasing the constant $\varepsilon_0$ in Lemma~\ref{lem:section-normalization} if necessary, there exist $\gamma_1\in(0,1)$ and $c,C>0$, depending only on $n,q,\lambda$, with the following property. Suppose that $v(b)=\min v$ and $a_0(v)\leq\varepsilon_0$. Then, for every $0<s<1$ satisfying $s^{-1}a_0(v)\leq\varepsilon_0$,
\[
b+c s^{\frac{1}{1+\gamma_1}}(K_1-b)\subset K_s\subset b+C s^{\frac{\gamma_1}{1+\gamma_1}}(K_1-b).
\]
If $v$ is normalized, then
\begin{equation}\label{eq:normalized-annular-growth}
c\abs{x}^{\frac{1+\gamma_1}{\gamma_1}}\leq v(x)\leq C\abs{x}^{1+\gamma_1}\quad\text{in }\left\{\varepsilon_0^{-1}a_0(v)<v<1\right\}.
\end{equation}
\end{Lemma}

\begin{proof}
Translate coordinates so that $b=0$. By Lemma~\ref{lem:section-normalization}, it suffices to consider normalized solutions. The Alexandrov estimate and the compactness argument in Lemma~\ref{lem:positive-level-geometry} give some $\varepsilon>0$ such that
\[
\frac{1+\varepsilon}{2}K_1\subset K_{1/2}\subset(1-\varepsilon)K_1.
\]
For every $0<s<1$ satisfying $s^{-1}a_0(v)\leq\varepsilon_0$, one has $a_0(v_{s,T})=s^{-1}a_0(v)\leq\varepsilon_0$ for every $T=T_{A,b}$ with $v(b)=\min v$. After normalization at successive heights, iterating the preceding inclusions as in the proof of \cite[Theorem~1.7]{JTX2025} yields the stated bounds for $K_s$. Together with $B_c\subset K_1\subset B_C$, these inclusions yield \eqref{eq:normalized-annular-growth}.
\end{proof}

\begin{Proposition}
\label{prop:normalized-annular-compactness}
Let $\{w_j\}$ be a sequence of normalized solutions on convex domains $\Omega_j$ exhausting $\R^n$, with $a_0(w_j)\to0$. Suppose that $w_j$ solves \eqref{eq:local-composite-equation} with coefficient $g_j$ satisfying, for some $\gamma_0\in(0,1)$,
\[
g_j\to1\quad\text{in }C^{\gamma_0}_{\mathrm{loc}}(\R^n).
\]
Then, after passing to a subsequence, $w_j$ converges in $C^1_{\mathrm{loc}}(\R^n)$ to a nonnegative convex solution
\[
v_\infty\in C^1_{\mathrm{loc}}(\R^n)\cap C^\infty_{\mathrm{loc}}(\R^n\setminus\{0\})
\]
of \eqref{eq:global-obstacle} on $\R^n$. There exist $\gamma_1\in(0,1)$ and $c,C>0$, depending only on $n,q,\lambda$, such that
\begin{equation}\label{eq:compactness-limit-growth}
c\min\left\{\abs{x}^{1+\gamma_1},\abs{x}^{\frac{1+\gamma_1}{\gamma_1}}\right\}\leq v_\infty(x)\leq C\max\left\{\abs{x}^{1+\gamma_1},\abs{x}^{\frac{1+\gamma_1}{\gamma_1}}\right\}.
\end{equation}
In particular, $\{v_\infty=0\}=\{0\}$. The convergence also holds in $C^{2,\gamma_0}_{\mathrm{loc}}(\R^n\setminus\{0\})$. If $g_j\equiv1$, then it holds in $C^\infty_{\mathrm{loc}}(\R^n\setminus\{0\})$.
\end{Proposition}

\begin{proof}
Set $\sigma_j=\varepsilon_0^{-1}a_0(w_j)$. Applying Lemma~\ref{lem:small-a0-annular-separation} to $w_j$ below unit height and to its normalized rescalings above unit height gives, for every $R>0$ and all sufficiently large $j$,
\[
c\min\left\{\abs{x}^{1+\gamma_1},\abs{x}^{\frac{1+\gamma_1}{\gamma_1}}\right\}\leq w_j(x)\leq C\max\left\{\abs{x}^{1+\gamma_1},\abs{x}^{\frac{1+\gamma_1}{\gamma_1}}\right\},\qquad x\in B_R\setminus K_{\sigma_j}(w_j).
\]
Since the $w_j$ are convex and $\sigma_j\to0$, these bounds give uniform local Lipschitz estimates. The Arzel\`a--Ascoli theorem yields, after passing to a subsequence, locally uniform convergence to a convex function $v_\infty$. Passing to the limit shows that $v_\infty$ solves \eqref{eq:global-obstacle} and satisfies \eqref{eq:compactness-limit-growth}. Lemma~\ref{lem:coincidence-set-regularity} gives $v_\infty\in C^1_{\mathrm{loc}}(\R^n)$, and convergence of gradients for convex functions yields $w_j\to v_\infty$ in $C^1_{\mathrm{loc}}(\R^n)$. Standard interior regularity gives $v_\infty\in C^\infty_{\mathrm{loc}}(\R^n\setminus\{0\})$ and the stated convergence away from the origin. Since the normalized condition gives $w_j(0)=\min_{\Omega_j}w_j$ and $0\leq w_j(0)\leq a_0(w_j)$, we have $v_\infty(0)=0$. The lower bound in \eqref{eq:compactness-limit-growth} is positive for $x\neq0$, and hence $\{v_\infty=0\}=\{0\}$.
\end{proof}

\subsection{Proof of Theorem~\ref{thm:complete-classification}}
\label{sec:proof-global-classification}

We prove Theorem~\ref{thm:complete-classification} by comparing the upper bound for $Q_n$ with its limit at infinity.

We first describe the radial solutions of \eqref{eq:entire-solution} for $0< q<n$. For a nontrivial nonnegative convex radial solution $w\in C^1(\R^n)$, we write $w(x)=\phi(|x|)$. The equation then gives
\begin{equation*}
\phi''(r)\left(\frac{\phi'(r)}{r}\right)^{n-1}=\phi(r)^q\quad\text{in }\{\phi>0\}.
\end{equation*}
Integrating, we obtain
\begin{equation*}
\phi'(r)^n=n\int_0^r \tau^{n-1}\phi(\tau)^q\,\ud\tau\quad \mbox{for all }r\geq0.
\end{equation*}
For $0<q<n$, standard ODE arguments give a unique positive entire profile $\phi_-$ with $\phi_-(0)=1$ and a unique entire profile $\phi_+$ vanishing on $[0,1]$ and positive on $(1,\infty)$. Comparison in the radial integral identity shows that the only profile vanishing precisely at the origin is $\phi(r)=c_*r^\frac{2n}{n-q}$. See also \cite{JiangTrudingerXu2026} for related ODE arguments and criteria for entire existence. For $q=0$, the same conclusion follows by direct integration. Scaling these profiles and including the homogeneous solution, we define
\[
\phi_a(r)=
\begin{cases}
a^{\frac{2n}{n-q}}\phi_+(r/a),&a>0,\\[3pt]
c_*r^{\frac{2n}{n-q}},&a=0,\\[3pt]
(-a)^{\frac{2n}{n-q}}\phi_-(r/(-a)),&a<0,
\end{cases}
\]
and set
\begin{equation}\label{eq:radial-family}
\Psi_a(x)=\phi_a(\abs{x}),\qquad x\in\R^n.
\end{equation}
For $a>0$, the coincidence set of $\Psi_a$ is $\overline{B_a}$; for $a<0$, $\Psi_a$ has minimum $(-a)^{\frac{2n}{n-q}}$ at $0$.

\begin{proof}[\textbf{Proof of Theorem~\ref{thm:complete-classification}}]
We first note that every nonempty sublevel set $K_s$ is bounded. Otherwise, since $K$ is bounded, after rotation and translation, we may assume that $K\cap\R^n_+=\emptyset$ and $K_s\cap\R^n_+$ is unbounded, where $\R^n_+=\{x_n>0\}$. By \cite[Lemma~2.5]{JTX2025}, the volume of $K_s\cap\R^n_+$ is bounded by $Cs^{(n-q)/2}$, whereas an unbounded open convex set has infinite volume, a contradiction.

Thus all sublevel sets are bounded, and $v$ attains its minimum. Lemma~\ref{lem:section-normalization} gives normalized rescalings $v_{s,T_s}$ as $s\to\infty$. Lemma~\ref{lem:positive-level-geometry} and \eqref{eq:Qn-scaling} give, for all sufficiently large $s$,
\[
Q_n[v](s)= Q_n[v_{2s,T_{2s}}](1/2) \in [c,C].
\]
By Theorem~\ref{thm:main-monotonicity}, $Q_n[v](s)$ is nonincreasing in $s$ and hence has the finite limit
\[
Q_{\infty}:=\lim_{s\to\infty}Q_n[v](s)\in [c,C].
\]
Proposition~\ref{prop:normalized-annular-compactness} gives a subsequence $v_{s_j,T_j}$ converging smoothly away from the origin to an entire solution $v_\infty$ of \eqref{eq:entire-solution}, with coincidence set $\{0\}$. By \eqref{eq:Qn-scaling},
\[
Q_n[v_\infty](r)=\lim_{j\to\infty}Q_n[v_{s_j,T_j}](r)=\lim_{j\to\infty}Q_n[v](s_jr)=Q_{\infty},\qquad \forall r>0.
\]
Thus $Q_n[v_\infty]$ is constant, and the equality case of Theorem~\ref{thm:main-monotonicity} implies $Q_{\infty}=c_n$. On the other hand, Lemma~\ref{lem:Qn-upper-bound} gives $Q_n[v](s)\leq c_n$ for $s>\min v$. Since $Q_n[v](s)$ is nonincreasing in $s$, we obtain $Q_n[v]\equiv c_n$ on $(\min v,\infty)$. The equality case of Theorem~\ref{thm:main-monotonicity} then shows that $v(A^{-1}x+b)$ is radial for some $A\in SL(n)$ and $b\in\R^n$. Since $v\in C^1$, uniqueness of the radial profiles gives
\[
v(x)=\Psi_a\bigl(A(x-b)\bigr)
\]
with $a=\left(\frac{\abs K}{\omega_n}\right)^{1/n} -(\min v)^{(n-q)/(2n)}$.
\end{proof}

\begin{Remark}\label{rem:rough-rhs}
Nearly the same argument also gives an analogue of the rigidity result of Mooney--Sen \cite[Theorem~1.2]{MooneySen2026} in higher dimensions. Let $n\geq2$, $0\leq q<n$, and let $v$ be a nonnegative convex solution of
\[
\det D^2v=f(v)\quad\text{in }\R^n,\qquad\lambda t^q\leq f(t)\leq\Lambda t^q\quad\text{for }t>0,
\]
where $f$ is a nonnegative Borel function and $0<\lambda\leq\Lambda$. If the coincidence set $\{v=0\}$ is bounded, then there exist $A\in SL(n)$ and $b\in\R^n$ such that $v(A^{-1}x+b)$ is radial. Indeed, the normalized blow-downs satisfy $\det D^2v_j=f_j(v_j)$, where $f_j(t):=s_j^{-q}f(s_jt)\rightharpoonup^\ast f_\infty(t)$ in $L^\infty_{\mathrm{loc}}((0,\infty))$ after passing to a subsequence. Using their local $C^1$ convergence away from $\{0\}$ and the coarea formula, one obtains $\det D^2v_\infty=f_\infty(v_\infty)$, and the rest of the argument is the same.
\end{Remark}

\section{Asymptotics near an isolated zero}
\label{sec:blow-up-normalization}

Let $0<q<n$, let $\Omega\subset\R^n$ be a convex domain, and let $v$ be a  convex solution of \eqref{eq:local-composite-equation}. Throughout this section, we assume that
\[
\{v=0\}=\{0\},\qquad K_1\Subset\Omega,\qquad g(0)=1,\qquad g\in C^{\gamma_0}(K_1)
\]
for some $\gamma_0\in(0,1)$. We consider $s\to0$, or equivalently $t=-\log s\to+\infty$.

For each sufficiently small $s>0$, choose $A_s\in SL(n)$ so that $v_{s,A_s}$ is normalized. Lemma~\ref{lem:small-a0-annular-separation} implies $K_s\subset B_{Cs^{\gamma_1/(1+\gamma_1)}}$ for small $s$. Since $B_c\subset s^{-\alpha}A_sK_s$, we also have $\|s^\alpha A_s^{-1}\|\leq C\operatorname{diam}K_s\to0$ as $s\to0$, where $\alpha=\frac{n-q}{2n}$. Thus the rescaled domains $\Omega_{s,A_s}$ exhaust $\R^n$ as $s\to0$, and, for each fixed $R$,
\[
\|g_{s,A_s}-1\|_{C^{\gamma_0}(B_R)}\leq C_Rs^{\gamma_g},\qquad \gamma_g=\frac{\gamma_0\gamma_1}{1+\gamma_1}>0.
\]
Proposition~\ref{prop:normalized-annular-compactness} and Theorem~\ref{thm:complete-classification} therefore show that every subsequential limit of these normalized rescalings as $s\to0$ is a homogeneous ellipsoidal solution.

Following the two-dimensional analysis of Choi--Choi--Kim \cite{CCK2024}, we study the system satisfied by the support functions of the level sets of $v$. For $0<s<1$, set $\mathbf h(\theta,s)=\mathbf h_{K_s}(\theta)$. On every compact positive level interval, interior regularity gives $\mathbf h>0$, $\mathbf h_s>0$, and $\nabla_{\Sn}^2\mathbf h+\mathbf hI>0$. As in the proof of Theorem~\ref{thm:main-monotonicity}, parameterizing $\partial K_s$ by its height and outer unit normal gives $Dv(X(\theta,s))=\theta/\mathbf h_s(\theta,s)$ and
\begin{equation}\label{eq:support-equation}
\mathbf{h}_{ss}+g(\theta,s)s^q\mathbf{h}_s^{n+2}\det\bigl(\nabla_{\Sn}^2\mathbf{h}+\mathbf{h}I\bigr)=0,
\end{equation}
where $g(\theta,s)=g(X(\theta,s))$; see \cite[Theorem~5.1]{LonginettiSalani2007} for the computation. For the homogeneous model $c_*|x|^\frac{2n}{n-q}$, the support function is $c_*^{-\alpha}s^\alpha$. We remove this scaling from $\mathbf h$ and $\mathbf h_s$ by setting
\begin{equation*}
t=-\log s,\qquad H(\theta,t)=\alpha s^{-\alpha}\mathbf{h}(\theta,s),\qquad Z(\theta,t)=s^{1-\alpha}\mathbf{h}_s(\theta,s).
\end{equation*}
Here $t$ increases as the height $s$ decreases. For $g\equiv1$, the homogeneous solution $c_*|x|^\frac{2n}{n-q}$ corresponds to the stationary pair $(h_0,h_0)$, where
\begin{equation*}
h_0=\alpha c_*^{-\alpha}=[\alpha^{n-1}(1-\alpha)]^{1/(2n)}.
\end{equation*}
Substituting the rescaled variables into \eqref{eq:support-equation} yields the nonautonomous system
\begin{equation}\label{eq:support-system}
\left\{
\begin{aligned}
H_t&=\alpha(H-Z), \\
Z_t&=-(1-\alpha)Z+\alpha^{-(n-1)}g\det\bigl(\nabla_{\Sn}^2H+HI\bigr)Z^{n+2}.
\end{aligned}
\right.
\end{equation}
Here
\[
g(\theta,t)=g\left(\frac{s^\alpha}{\alpha}(H\theta+\nabla_{\Sn}H)\right).
\]
For each $t$, $H(\cdot,t)$ is the support function of $K_H:=\alpha s^{-\alpha}\overline{K_s}$.

For $s>0$ and $T=T_{A,b}$, with $A\in SL(n)$ and $b\in\R^n$, let $(H_{s,T},Z_{s,T})$ be the pair associated with $v_{s,T}$. Then
\begin{equation*}
(H_{s,T},Z_{s,T})(\cdot,\tau)=A[(H,Z)](\cdot,t')-\bigl(\alpha e^{\alpha t'}(Ab)\cdot\theta,0\bigr),
\end{equation*}
where $t'=\tau-\log s$ and, for functions $w$ on $\Sn$, we define
\[
A[w](\theta):=\abs{A^\top\theta}\,w\!\left(\frac{A^\top\theta}{\abs{A^\top\theta}}\right).
\]
The action on pairs is understood componentwise, and $A[\mathbf h_K]=\mathbf h_{AK}$. The transformed pair satisfies \eqref{eq:support-system} with $g_{s,T}(y)=g(b+s^\alpha A^{-1}y)$ in place of $g$. For $s=1$ and $b=0$, we write $g_A=g\circ A^{-1}$.

By the interior estimates on annuli between fixed positive levels in Lemma~\ref{lem:positive-level-geometry}, the convergence above gives matrices $\bar A(t)\in SL(n)$ such that
\begin{equation}\label{eq:local-subsequential-pair-convergence}
\bar A(t)[(H,Z)](\cdot,t+\cdot)\longrightarrow(h_0,h_0)\qquad\text{as }t\to+\infty
\end{equation}
in $C^{2,\gamma_0}(\Sn\times[-\tau_0,\tau_0])\times C^{1,\gamma_0}(\Sn\times[-\tau_0,\tau_0])$ for every fixed $\tau_0>0$. Let $\sigma(t)$ be the sum of the $C^{2,\gamma_0}$ and $C^{1,\gamma_0}$ norms of the differences from $(h_0,h_0)$ on $\Sn\times[-2,2]$; then $\sigma(t)\to0$ as $t\to+\infty$.

For each $\ell\geq0$, let $\mathcal H_\ell$ be the space of spherical harmonics of degree $\ell$, and let $P_\ell$ be the $L^2(\Sn)$-orthogonal projection onto $\mathcal H_\ell$. Then, for $m=0,1$, every $w\in H^m(\Sn)$ has the orthogonal expansion
\[
w=\sum_{\ell=0}^{\infty}w_\ell=\sum_{\ell=0}^{\infty}P_\ell w\qquad\text{with convergence in }H^m(\Sn),
\]
where $w_\ell=P_\ell w$ and $H^0(\Sn)=L^2(\Sn)$. The eigenvalue of $-\Delta_{\Sn}$ on $\mathcal H_\ell$ is $\lambda_\ell=\ell(\ell+n-2)$. Set $P_{\geq3}=I-P_0-P_1-P_2$ and, for $m\geq0$, write $H^m_{\geq3}(\Sn)=P_{\geq3}H^m(\Sn)$.

\subsection{Affine normalization}
\label{sec:normalization-spectrum}

Let $\operatorname{Sym}_0(n)$ denote the space of symmetric trace-free matrices, and set $q_B(\theta)=\theta^\top B\theta$ for $B\in\operatorname{Sym}_0(n)$. Then
\begin{equation*}
\mathcal H_2=\{q_B:B\in\operatorname{Sym}_0(n)\}.
\end{equation*}
The affine action gives
\[
\left.\frac{\ud}{\ud\varepsilon}\right|_{\varepsilon=0}e^{\varepsilon B}[(h_0,h_0)]=h_0(q_B,q_B).
\]
Thus the pairs $(Y_2,Y_2)$, $Y_2\in\mathcal H_2$, form the tangent space at $(h_0,h_0)$ to its $SL(n)$-orbit, which we call the ellipsoid orbit.

To choose these normalizations consistently, we place the rescaled level sets in minimal mean width position. For a convex body $K\subset\R^n$, define its mean width by
\begin{equation*}
w(K):=\frac{2}{\abs{\Sn}}\int_{\Sn}\mathbf h_K(\theta)\,\ud\theta.
\end{equation*}
We say that $K$ is in minimal mean width position if
\begin{equation*}
w(K)\leq w(AK)\qquad\text{for every }A\in SL(n).
\end{equation*}
By Giannopoulos and Milman \cite[Theorems~3.1 and~3.3]{GiannopoulosMilman2000}, this is equivalent to
\begin{equation*}
\int_{\Sn}\mathbf h_K(\theta)\,\theta\otimes\theta\,\ud\theta=\frac1n\left(\int_{\Sn}\mathbf h_K(\theta)\,\ud\theta\right)I,
\end{equation*}
and the choice of $A\in SL(n)$ that places $K$ in minimal mean width position is unique up to left multiplication by a rotation. This isotropy condition is equivalent to $P_2\mathbf h_K=0$. Accordingly, we call $A\in SL(n)$ a degree-two normalizing matrix for the support function $\mathbf h=\mathbf h_K$ of a convex body $K$ if $AK$ is in minimal mean width position. For $K=K_s$ and $t=-\log s$, this means
\[
P_2\left(A[\mathbf h]-h_0\right)=0,\quad\text{or equivalently}\quad P_2\{A[H]-h_0\}=0.
\]

\begin{Lemma}\label{lem:local-compatible-normalization}
There exist $t_1<\infty$ and a $C^1$ family $A:[t_1,\infty)\to SL(n)$ such that, with
\begin{equation}\label{eq:normalized-perturbations}
A(t)[H]=h_0(1+u),\qquad A(t)[Z]=h_0(1+z),
\end{equation}
we have $P_2u=0$, $L(t):=A_tA^{-1}\in\operatorname{Sym}_0(n)$, and
\begin{equation}\label{eq:normalization-space-time-convergence}
\|u\|_{C^{2,\gamma_0}([t-1,t+1]\times\Sn)}+\|z\|_{C^{1,\gamma_0}([t-1,t+1]\times\Sn)}\longrightarrow0\qquad\text{as }t\to+\infty.
\end{equation}
\end{Lemma}

\begin{proof}
For $S\in\operatorname{Sym}_0(n)$, set $\Psi(S,F)=P_2\{e^S[F]-h_0\}$. Since $D_S\Psi(0,h_0)[B]=h_0q_B$ is an isomorphism from $\operatorname{Sym}_0(n)$ onto $\mathcal H_2$, the implicit function theorem gives a unique small symmetric correction $S=S(F)$ for $F$ near $h_0$ in $C^{2,\gamma_0}(\Sn)$. Together with \eqref{eq:local-subsequential-pair-convergence}, this gives local $C^1$ families of degree-two normalizing matrices on overlapping intervals of fixed length for all sufficiently large $t$.

For such a local family $\bar A(t)$, set $M=\bar A_t\bar A^{-1}$ and solve
\[
R_t=-R\operatorname{skew}M
\]
with prescribed initial data in $SO(n)$. Then $A=R\bar A$ satisfies
\[
A_tA^{-1}=R\operatorname{sym}(M)R^{-1}\in\operatorname{Sym}_0(n).
\]
Rotations preserve the degree-two normalization. On an overlap, uniqueness of minimal mean width position gives $A_2=RA_1$ with $R$ taking values in $SO(n)$, and
\[
L_2=R_tR^{-1}+RL_1R^{-1},\qquad L_i=(A_i)_tA_i^{-1}.
\]
Since both $L_i$ are symmetric, $R_t=0$. The local families therefore agree after constant left rotations and extend to $[t_1,\infty)$.

Fix $t$ sufficiently large, choose a constant matrix $\bar A$ on $[t-2,t+2]$ as in \eqref{eq:local-subsequential-pair-convergence}, and set $F(\theta,\tau)=\bar A[H](\theta,t+\tau)$ for $|\tau|\leq2$. A change of variables on $\Sn$ shows that $\Psi$ is smooth also for $F\in C^0(\Sn)$. Its derivative in $S$ at $(0,h_0)$ is invertible. The implicit function theorem and its parameter estimates give
\[
\|S\|_{C^{2,\gamma_0}([-2,2])}\leq C\sigma(t).
\]
Set $M=(e^S)_\tau e^{-S}$ and solve $R_\tau=-R\operatorname{skew}M$ with $R(0)=I$. Since $\|M\|_{C^{1,\gamma_0}([-2,2])}\leq C\sigma(t)$, the rotation equation gives
\[
\|R-I\|_{C^{2,\gamma_0}([-2,2])}\leq C\sigma(t).
\]
Consequently, the local family $R(\tau)e^{S(\tau)}\bar A$ preserves the $C^{2,\gamma_0}$ convergence of $H$ and the $C^{1,\gamma_0}$ convergence of $Z$ on the space-time cylinder. The global family differs from this local family by a constant left rotation and therefore satisfies the same convergence estimates. This proves \eqref{eq:normalization-space-time-convergence}. The corresponding rescalings converge locally uniformly to the round homogeneous solution as $s=e^{-t}\to0$.
\end{proof}

\subsection{The normalized system and spectral analysis}

For a smooth family $w(\cdot,t)$ on $\Sn$, differentiation under the $SL(n)$ action gives
\begin{equation*}
\partial_t(A[w])=A[w_t]+\Phi_L[A[w]],\qquad L=A_tA^{-1}.
\end{equation*}
Here, for a trace-free matrix $B$ and $w\in C^1(\Sn)$, we define
\[
\Phi_B[w](\theta):=\left.\frac{\ud}{\ud\varepsilon}\right|_{\varepsilon=0}e^{\varepsilon B}[w](\theta)=\left\langle B\big(w(\theta)\theta+\nabla_{\Sn}w(\theta)\big),\theta\right\rangle .
\]
For $L\in\operatorname{Sym}_0(n)$, we have $\Phi_L[h_0]=h_0q_L$. More generally, for a scalar $c$ and a smooth function $w$ on $\Sn$,
\[
\Phi_L[c]=c q_L,\qquad\|\Phi_L[w]\|_{H^m}\leq C_{m,n}|L|\|w\|_{H^{m+1}}\qquad(m\geq0).
\]
Differentiating \eqref{eq:normalized-perturbations} and using \eqref{eq:support-system}, we obtain
\begin{equation}\label{eq:normalized-system}
\left\{
\begin{aligned}
u_t&=\alpha(u-z)+q_L+\Phi_L[u],\\
z_t&=(1-\alpha)\bigl[-(1+z)+\widetilde g\,\det\bigl(I+\nabla_{\Sn}^2u+uI\bigr)(1+z)^{n+2}\bigr]\\
&\quad+q_L+\Phi_L[z],
\end{aligned}
\right.
\end{equation}
where $\widetilde g(\theta,t)=g_{A(t)}(\theta,t)$.

Lemma~\ref{lem:small-a0-annular-separation} gives $K_s\subset B_{Cs^{\gamma_1/(1+\gamma_1)}}$ for small $s>0$, so the H\"older continuity of $g$ and $g(0)=1$ yield
\begin{equation}\label{eq:normalized-coefficient-decay}
\|\widetilde g(\cdot,t)-1\|_{L^\infty(\Sn)}\leq Ce^{-\gamma_gt}\qquad\text{as }t\to+\infty,
\end{equation}
with $\gamma_g=\gamma_0\gamma_1/(1+\gamma_1)$.

Set $z_2=P_2z$. Since $P_2u_t=P_2u=0$, projecting the first equation in \eqref{eq:normalized-system} onto $\mathcal H_2$ gives
\begin{equation}\label{eq:degree-two-modulation}
q_L+P_2\Phi_L[u]=\alpha z_2.
\end{equation}
It follows from \eqref{eq:normalization-space-time-convergence} that, after increasing $t_1$ if necessary, the map
\[
L\longmapsto q_L+P_2\Phi_L[u]
\]
is an isomorphism from $\operatorname{Sym}_0(n)$ onto $\mathcal H_2$, with uniformly bounded inverse. Hence
\begin{equation}\label{eq:local-affine-velocity}
|L(t)|\leq C(n,q)\|P_2z(t)\|_{L^2}\qquad(t\geq t_1).
\end{equation}
Equation~\eqref{eq:degree-two-modulation} uniquely determines $L$, which depends smoothly on $(u,z_2)$ and linearly on $z_2$ for fixed $u$. Thus an integrable bound for $z$ will give convergence of the normalizing matrices.

We now separate the terms containing $\widetilde g-1$ and expand the autonomous part of \eqref{eq:normalized-system} at $(u,z)=(0,0)$. Using the expansion
\[
\det\bigl(I+\nabla_{\Sn}^2u+uI\bigr)=1+\Delta_{\Sn}u+(n-1)u+O\left((\abs{u}+\abs{\nabla_{\Sn}^2u})^2\right),
\]
and the linear part $q_L=\alpha P_2z$ of \eqref{eq:degree-two-modulation}, we obtain the linear operator
\begin{equation}\label{eq:forward-linear-operator}
\mathbf A_{\rm red}\binom{u}{z}=\binom{\alpha(u-(I-P_2)z)}{(1-\alpha)\bigl((\Delta_{\Sn}+n-1)u+(n+1)z\bigr)+\alpha P_2z}
\end{equation}
on the domain $\bigl(\mathcal H_0\oplus\mathcal H_1\oplus H^2_{\geq3}(\Sn)\bigr) \times H^1(\Sn)$. Substituting \eqref{eq:degree-two-modulation} into \eqref{eq:normalized-system} gives
\begin{equation*}
(u_t,z_t)=\mathbf A_{\rm red}(u,z)+\mathbf N_{\rm red}(u,z)+\mathbf F_g(u,z,t),
\end{equation*}
where
\begin{equation*}
\mathbf N_{\rm red}(u,z):=\binom{(I-P_2)\Phi_L[u]}{(1-\alpha)\mathcal R_{\rm det}(u,z)+\Phi_L[z]-P_2\Phi_L[u]},
\end{equation*}
\begin{align*}
\mathcal R_{\rm det}(u,z)& =-(1+z)+\det\bigl(I+\nabla_{\Sn}^2u+uI\bigr)(1+z)^{n+2}\\
&\quad-(\Delta_{\Sn}+n-1)u-(n+1)z,
\end{align*}
and
\[
\mathbf F_g(u,z,t):=\binom{0}{(1-\alpha)(\widetilde g-1)\det\bigl(I+\nabla_{\Sn}^2u+uI\bigr)(1+z)^{n+2}}.
\]
The limiting linearized operator is therefore the same as for $g\equiv1$, and $\mathbf F_g$ contains all terms involving $\widetilde g-1$.

We next compute the spectrum of $\mathbf A_{\rm red}$. On spherical harmonics of degree $\ell\ne2$, the linearized equations are
\[
\left\{
\begin{aligned}
u_t&=\alpha(u-z), \\
z_t&=(1-\alpha)\bigl((\Delta_{\Sn}+n-1)u+(n+1)z\bigr).
\end{aligned}
\right.
\]
Eliminating $z$ gives
\[
u_{tt}-\kappa u_t+\alpha(1-\alpha)(\Delta_{\Sn}+2n)u=0,
\]
where
$$
\kappa=n+1-n\alpha=\frac{n+q+2}{2}.
$$
Projecting this equation onto $\mathcal H_\ell$ gives modes of the form $u(\theta,t)=e^{\mu t}Y_\ell(\theta)$ with $0\ne Y_\ell\in\mathcal H_\ell$ and $\mu=\mu_\ell^\pm$, where
\begin{equation}\label{eq:indicial-roots}
\mu_\ell^\pm=\frac{\kappa\pm\sqrt{\kappa^2+4\alpha(1-\alpha)(\lambda_\ell-2n)}}{2}.
\end{equation}
For $\ell=0,1$, the roots satisfy $0<\mu_\ell^-<\mu_\ell^+$. For $\ell=1$, the roots are $\alpha$ and $\kappa-\alpha$, and the modes $(u,z)=(e^{\alpha t}Y_1,0)$ are generated by translations.
The remaining roots satisfy
\[
\cdots<\mu_4^-<\mu_3^-<0<\kappa<\mu_3^+<\mu_4^+<\cdots.
\]

In degree two, the normalization $P_2u=0$ removes the zero modes, and the normalization terms in \eqref{eq:forward-linear-operator} give
\[
\mathbf A_{\rm red}\binom{0}{Y_2}=\kappa\binom{0}{Y_2},\qquad Y_2\in\mathcal H_2.
\]

Define
\[
\mathcal X_{\rm red}=\bigl(\mathcal H_0\oplus\mathcal H_1\oplus H^1_{\geq3}(\Sn)\bigr)\times L^2(\Sn).
\]
To obtain orthogonal spectral projections, we use the equivalent Hilbert norm
\[
\begin{aligned}
\|(u,z)\|_{\mathcal X_{\rm red}}^2={}&\frac{1-\alpha}{\alpha}\left(\|\nabla_{\Sn}u\|_{L^2}^2-(n-1)\|u\|_{L^2}^2\right)+\|z\|_{L^2}^2\\
&+\frac{\kappa-2\alpha}{\alpha}\sum_{\ell=0}^{1}\left[\langle P_\ell u,P_\ell z\rangle_{L^2}+\frac{\kappa-2\alpha}{2\alpha}\|P_\ell u\|_{L^2}^2\right],
\end{aligned}
\]
and let $\langle\cdot,\cdot\rangle_{\mathcal X_{\rm red}}$ denote the associated inner product.

\begin{Lemma}
\label{lem:forward-symmetrizer}
The norm $\|\cdot\|_{\mathcal X_{\rm red}}$ is equivalent to the standard $H^1\times L^2$ norm on $\mathcal X_{\rm red}$, with equivalence constants depending only on $n$ and $q$. With respect to the associated inner product, $\mathbf A_{\rm red}$ is self-adjoint on the stated domain. The space $\mathcal X_{\rm red}$ has the orthogonal spectral decomposition
\[
\mathcal X_{\rm red}=\mathcal X_+\oplus\mathcal X_-.
\]
Here $\mathcal X_+$ and $\mathcal X_-$ are the positive and negative spectral subspaces of $\mathbf A_{\rm red}$, respectively. Moreover, in the sense of quadratic forms,
\[
\mathbf A_{\rm red}|_{\mathcal X_+}\geq\alpha I,\qquad\mathbf A_{\rm red}|_{\mathcal X_-}\leq\mu_3^-I.
\]
\end{Lemma}

\begin{proof}
For $\ell\geq0$, $\ell\ne2$, set
\[
\mathcal E_{\ell,\pm}=\left\{\left(Y_\ell,\left(1-\frac{\mu_\ell^\pm}{\alpha}\right)Y_\ell\right):Y_\ell\in\mathcal H_\ell\right\},
\]
and set
\[
\mathcal E_{2,+}=\{(0,Y_2):Y_2\in\mathcal H_2\}.
\]
By \eqref{eq:forward-linear-operator}, $\mathbf A_{\rm red}$ acts on $\mathcal E_{\ell,\pm}$ by multiplication by $\mu_\ell^\pm$, and on $\mathcal E_{2,+}$ by multiplication by $\kappa$. Both roots are positive for $\ell=0,1$, whereas $\mu_\ell^-<0<\mu_\ell^+$ for $\ell\geq3$. Thus
\[
\mathcal X_+=(\mathcal E_{0,-}\oplus\mathcal E_{1,-})\oplus\overline{\bigoplus_{\ell=0}^{\infty}\mathcal E_{\ell,+}},\qquad\mathcal X_-=\overline{\bigoplus_{\ell=3}^{\infty}\mathcal E_{\ell,-}},
\]
where the closures are taken in $\mathcal X_{\rm red}$.

For $\ell=0,1$, the degree-$\ell$ contribution to $\|(u,z)\|_{\mathcal X_{\rm red}}^2$ is
\[
\left\|P_\ell z+\frac{\kappa-2\alpha}{2\alpha}P_\ell u\right\|_{L^2}^2+\frac{(\mu_\ell^+-\mu_\ell^-)^2}{4\alpha^2}\|P_\ell u\|_{L^2}^2.
\]
Since $\mu_\ell^+\ne\mu_\ell^-$, this quadratic form is positive definite. Substituting the eigenvectors above also shows that $\mathcal E_{\ell,+}$ and $\mathcal E_{\ell,-}$ are orthogonal. For $\ell\geq3$, set
\[
a_\ell=\frac{1-\alpha}{\alpha}\bigl(\lambda_\ell-(n-1)\bigr)>0.
\]
The degree-$\ell$ contribution to the squared norm is $a_\ell\|P_\ell u\|_{L^2}^2+\|P_\ell z\|_{L^2}^2$, and
\[
\left(1-\frac{\mu_\ell^+}{\alpha}\right)\left(1-\frac{\mu_\ell^-}{\alpha}\right)=-a_\ell.
\]
Hence $\mathcal E_{\ell,+}$ and $\mathcal E_{\ell,-}$ are orthogonal. The degree-two contribution is $\|P_2z\|_{L^2}^2$, and different spherical harmonic degrees are orthogonal. Since $a_\ell$ is comparable to $1+\lambda_\ell$ for $\ell\geq3$, these identities and the calculation for $\ell=0,1,2$ prove the norm equivalence and the completeness of the orthogonal eigenspace decomposition.

Write $U_{\ell,\pm}$ for the components of $U\in\mathcal X_{\rm red}$ in $\mathcal E_{\ell,\pm}$. The maximal domain of the resulting real diagonal operator is
\[
\left\{U\in\mathcal X_{\rm red}:\sum_{\ell\geq3}\sum_{\pm}|\mu_\ell^\pm|^2\|U_{\ell,\pm}\|_{\mathcal X_{\rm red}}^2<\infty\right\}.
\]
Since $|\mu_\ell^\pm|$ is comparable to $(1+\lambda_\ell)^{1/2}$ for $\ell\geq3$, this domain is precisely
\[
(\mathcal H_0\oplus\mathcal H_1\oplus H^2_{\geq3})\times H^1,
\]
and its graph norm is equivalent to the $H^2\times H^1$ norm. A real diagonal operator on its maximal domain is self-adjoint, proving the assertion for $\mathbf A_{\rm red}$. The spectral splitting and bounds follow from \eqref{eq:indicial-roots} and the strict decrease of $\mu_\ell^-$.
\end{proof}

Let $\Pi_+$ and $\Pi_-$ denote the orthogonal spectral projections of $\mathbf A_{\rm red}$ onto its positive and negative spectral subspaces $\mathcal X_+$ and $\mathcal X_-$. Write $U=(u,z)$ and set
\[
R(t)=\mathbf N_{\rm red}(u(t),z(t))+\mathbf F_g(u(t),z(t),t).
\]
The Taylor expansion of $\mathcal R_{\rm det}$, \eqref{eq:local-affine-velocity}, and \eqref{eq:normalization-space-time-convergence} give
\[
\begin{aligned}
\|\mathbf N_{\rm red}(u(t),z(t))\|_{\mathcal X_{\rm red}}&\leq \varepsilon(t)\bigl(\|u(t)\|_{H^2}+\|z(t)\|_{H^1}\bigr)\\
&\leq C\varepsilon(t)\bigl(\|\mathbf A_{\rm red}U(t)\|_{\mathcal X_{\rm red}}+\|U(t)\|_{\mathcal X_{\rm red}}\bigr),
\end{aligned}
\]
where $\varepsilon(t)\to0$ as $t\to+\infty$. The second inequality follows from the equivalence of the graph norm of $\mathbf A_{\rm red}$ and the $H^2\times H^1$ norm, proved in Lemma~\ref{lem:forward-symmetrizer}. Therefore, \eqref{eq:normalized-coefficient-decay} yields
\[
\|R(t)\|_{\mathcal X_{\rm red}}\leq \varepsilon(t)\bigl(\|\mathbf A_{\rm red}U(t)\|_{\mathcal X_{\rm red}}+\|U(t)\|_{\mathcal X_{\rm red}}\bigr)+Ce^{-\gamma_gt},
\]
after enlarging $\varepsilon$ by a fixed factor.

The estimate for $R$ involves $\mathbf A_{\rm red}U$, so we first bound the weighted $L^2$ norm in time of $\mathbf A_{\rm red}U$ by that of $U$. In the following estimates, norms and inner products are taken in $\mathcal X_{\rm red}$. Fix $0<\sigma<\min\{\gamma_g,-\mu_3^-\}$ and $a>2\gamma_g$, and take $T$ sufficiently large. Since $U_t-\mathbf A_{\rm red}U=R$, self-adjointness and integration by parts give
\[
\begin{aligned}
&\int_T^\infty e^{-a|t-s|}\bigl(\|U_t\|^2+\|\mathbf A_{\rm red}U\|^2\bigr)\,\ud s\\
&\quad=\int_T^\infty e^{-a|t-s|}\|R\|^2\,\ud s+2\int_T^\infty e^{-a|t-s|}\langle U_t,\mathbf A_{\rm red}U\rangle\,\ud s\\
&\quad=\int_T^\infty e^{-a|t-s|}\|R\|^2\,\ud s-\int_T^\infty \partial_s\!\left(e^{-a|t-s|}\right)\langle U,\mathbf A_{\rm red}U\rangle\,\ud s\\
&\qquad\quad-e^{-a(t-T)}\langle U(T),\mathbf A_{\rm red}U(T)\rangle.
\end{aligned}
\]
Using $\gamma_g>\sigma$ and $a>2\sigma$, we obtain, for all sufficiently large $t$,
\begin{equation}\label{eq:forward-weighted-graph-estimate}
\int_T^\infty e^{-a|t-s|}\bigl(\|U_t\|^2+\|\mathbf A_{\rm red}U\|^2\bigr)\,\ud s\leq C_a\left(\int_T^\infty e^{-a|t-s|}\|U\|^2\,\ud s+e^{-2\sigma t}\right).
\end{equation}

Define
\[
\begin{aligned}
N_+(t)^2&=\int_T^\infty e^{-a|t-s|}\|\Pi_+U\|^2\,\ud s,\\
N_-(t)^2&=\int_T^\infty e^{-a|t-s|}\|\Pi_-U\|^2\,\ud s+e^{-2\sigma t}.
\end{aligned}
\]
By orthogonality, the right-hand side of \eqref{eq:forward-weighted-graph-estimate} is $C_a(N_+^2+N_-^2)$, which tends to zero as $t\to+\infty$. The estimate on $R$ implies
\[
r(t):=\left(\int_T^\infty e^{-a|t-s|}\|R\|^2\,\ud s\right)^{1/2}=o\bigl(N_+(t)+N_-(t)\bigr)\qquad\text{as }t\to+\infty.
\]
Using the spectral bounds in Lemma~\ref{lem:forward-symmetrizer}, differentiation and integration by parts give
\[
\begin{aligned}
N_+'(t)&\geq\alpha N_+(t)-r(t),\\
N_-'(t)&\leq-\sigma N_-(t)+r(t)+C_Te^{-(a-\sigma)t}.
\end{aligned}
\]
The last term is $o(N_-)$ since $a>2\sigma$.

We use the following form of the Merle--Zaag lemma \cite{MZ1998}; see \cite[Lemma~A.1]{CCK2024}.

\begin{Lemma}[Merle--Zaag lemma]
\label{lem:three-component-MZ}
Let $x,y,z$ be nonnegative locally absolutely continuous functions on $[T,\infty)$ such that
\[
x(t)+y(t)+z(t)>0,\qquad x(t)+y(t)+z(t)\rightarrow0\qquad\text{as }t\to+\infty.
\]
Suppose that there exists $c_0>0$ such that, for every $\varepsilon>0$, there exists $T_\varepsilon\geq T$ for which the following inequalities hold for almost every $t\geq T_\varepsilon$:
\[
\left\{
\begin{aligned}
x'-c_0x&\geq-\varepsilon(y+z),\\
\abs{y'}&\leq\varepsilon(x+y+z),\\
z'+c_0z&\leq\varepsilon(x+y).
\end{aligned}
\right.
\]
Then, as $t\to+\infty$, one of the following alternatives holds:
\[
x+z=o(y)\qquad\text{or}\qquad x+y=o(z).
\]
\end{Lemma}

\begin{Proposition}\label{thm:forced-spectral-decay}
For every $0<\mu<\min\{\gamma_g,-\mu_3^-\}$,
\begin{equation}\label{eq:forced-U-decay}
\|(u(t),z(t))\|_{\mathcal X_{\rm red}}+|L(t)|\leq C_\mu e^{-\mu t}\qquad\text{as }t\to+\infty.
\end{equation}
Consequently, $A(t)\to A_\infty$ in $SL(n)$ as $t\to+\infty$.
\end{Proposition}

\begin{proof}
Fix $0<\mu<\min\{\gamma_g,-\mu_3^-\}$. Choose $\sigma$ with $\mu<\sigma<\min\{\gamma_g,-\mu_3^-\}$ and $a>2\gamma_g$, and let $N_\pm$ be defined as above. The preceding differential inequalities and Lemma~\ref{lem:three-component-MZ}, applied to $(N_+,0,N_-)$, give $N_+=o(N_-)$ as $t\to+\infty$. Substituting this relation into the inequality for $N_-$ and integrating, we obtain
\[
N_-'\leq-(\sigma-o(1))N_-,\qquad N_+(t)+N_-(t)\leq C_\mu e^{-\mu t}
\]
for all sufficiently large $t$.

The Sobolev embedding in the time variable and \eqref{eq:forward-weighted-graph-estimate} now give
\[
\|U(t)\|^2\leq C\int_{t-1}^{t+1}\bigl(\|U\|^2+\|U_t\|^2\bigr)\,\ud s\leq C_a\bigl(N_+(t)^2+N_-(t)^2\bigr)\leq C_\mu e^{-2\mu t}.
\]
The preceding estimate, together with $|L|\leq C\|z_2\|_{L^2}\leq C\|(u,z)\|_{\mathcal X_{\rm red}}$, proves \eqref{eq:forced-U-decay}. In particular, $L$ is integrable. Gronwall's inequality applied to $A_t=LA$ gives a uniform bound for $A$, and integration then yields $A(t)\to A_\infty\in SL(n)$ as $t\to+\infty$.
\end{proof}

\subsection{Asymptotics in fixed affine coordinates}

\begin{proof}[\textbf{Proof of Theorem~\ref{thm:intro-local-leading}}]
Normalizing at a sufficiently small height reduces the problem to the setting of this section. The asserted asymptotic estimate is preserved by this rescaling. By \eqref{eq:normalization-space-time-convergence} and Proposition~\ref{thm:forced-spectral-decay},
\[
A(t)\longrightarrow A_\infty\in SL(n),\qquad A(t)[H](\cdot,t)\longrightarrow h_0\quad\text{as }t\to+\infty.
\]
Since the support functions converge uniformly on $\Sn$ and the $SL(n)$ action is continuous in $C^0(\Sn)$, the convergence of $A(t)$ gives
\[
A_\infty[H](\cdot,t)=h_0+o(1)\qquad\text{as }t\to+\infty.
\]
Let $w(x)=v(A_\infty^{-1}x)$. It follows that the radial function of $K_s(w)$ is
\[
\frac{h_0}{\alpha}s^\alpha(1+o(1))\qquad\text{as }s\to0,
\]
uniformly in the angular variable. Inverting this relation gives
\[
w(x)=\left(\frac{\alpha}{h_0}\right)^\frac{2n}{n-q}\abs{x}^\frac{2n}{n-q}+o(\abs{x}^\frac{2n}{n-q})=c_*\abs{x}^\frac{2n}{n-q}+o(\abs{x}^\frac{2n}{n-q}).
\]
The preceding asymptotic, convexity, and the H\"older continuity of $g$ give $g(x)-1=O(\abs{x}^\gamma)$ as $x\to0$, and hence $\widetilde g-1=O(e^{-\alpha\gamma t})$ as $t\to+\infty$. With this improved coefficient decay, \eqref{eq:forced-U-decay} holds with $\gamma_g$ replaced by $\alpha\gamma$. Since $\alpha\gamma<-\mu_3^-$, standard interior elliptic estimates in fixed affine coordinates then yield
\[
\|u(t)\|_{C^2}+\|z(t)\|_{C^1}=O(e^{-\mu t})
\]
as $t\to+\infty$ for every $0<\mu<\alpha\gamma$.

Choose $\alpha\gamma/2<\mu<\alpha\gamma$. The quadratic structure of $\mathbf N_{\rm red}$ gives
\[
\|\mathbf N_{\rm red}(u,z)+\mathbf F_g(u,z,t)\|_{\mathcal X_{\rm red}}\leq C\bigl(e^{-2\mu t}+e^{-\alpha\gamma t}\bigr)\leq Ce^{-\alpha\gamma t}.
\]
Integrating the projected equations and using the bound for $L$ in $A_t=LA$, we obtain
\[
\|(u(t),z(t))\|_{\mathcal X_{\rm red}}+|A(t)-A_\infty|=O(e^{-\alpha\gamma t})\qquad\text{as }t\to+\infty.
\]
Applying the same interior estimates and inverting the level-set relation again gives
\[
w(x)=c_*\abs{x}^\frac{2n}{n-q}+O(\abs{x}^{\frac{2n}{n-q}+\gamma})\qquad\text{as }x\to0.
\]
This proves \eqref{eq:intro-v-leading}.
\end{proof}

\noindent \textbf{AI declaration.} The proof strategy using a monotonicity formula, level-set parametrization, and the Merle--Zaag lemma was developed by the authors. We used GPT-6 Astra to assist with some technical derivations in convex geometry and spectral analysis. The proof of Theorem~\ref{thm:complete-classification} presented here builds on the monotone quantity \eqref{eq:Qn-introduction}, which is a natural extension to higher dimensions of the quantity obtained in dimension two in Mooney--Sen \cite{MooneySen2026}.  We wrote the manuscript after working out the proofs and used GPT-6 Astra to help polish the language.

\bigskip
\noindent
Tianling Jin\\
Department of Mathematics, The Hong Kong University of Science and Technology\\
Clear Water Bay, Kowloon, Hong Kong\\
\textit{Email address:} \texttt{tianlingjin@ust.hk}

\medskip
\noindent
Connor Mooney\\
Department of Mathematics, UC Irvine\\
\textit{Email address:} \texttt{mooneycr@uci.edu}

\medskip
\noindent
Aranya Sen\\
Department of Mathematics, UC Irvine\\
\textit{Email address:} \texttt{aranyas1@uci.edu}

\medskip
\noindent
Xushan Tu\\
Department of Mathematics, The Hong Kong University of Science and Technology\\
Clear Water Bay, Kowloon, Hong Kong\\
\textit{Email addresses:} \texttt{maxstu@ust.hk}, \texttt{tuxushan28561m@gmail.com}

\medskip
\noindent
Jingang Xiong\\
School of Mathematical Sciences, Laboratory of Mathematics and Complex Systems, MOE\\
Beijing Normal University, Beijing 100875, China\\
\textit{Email address:} \texttt{jx@bnu.edu.cn}


\begin{thebibliography}{99}
\small
\setlength{\itemsep}{0pt}

\bibitem{AltPhillips1986}
H.~W.~Alt and D.~Phillips,
\newblock \emph{A free boundary problem for semilinear elliptic equations},
\newblock J. Reine Angew. Math. \textbf{368} (1986), 63--107.

\bibitem{bonorino2001regularity}
L.~P.~Bonorino,
\newblock \emph{Regularity of the free boundary for some elliptic and parabolic problems. I},
\newblock Comm. Partial Differential Equations \textbf{26} (2001), no.~1--2, 175--203.

\bibitem{BonorinoEtAl2005}
L.~P.~Bonorino, E.~H.~M.~Brietzke, J.~P.~Lukaszczyk, and C.~A.~Taschetto,
\newblock \emph{Properties of the period function for some Hamiltonian systems and homogeneous solutions of a semilinear elliptic equation},
\newblock J. Differential Equations \textbf{214} (2005), no.~1, 156--175.

\bibitem{Caffarelli1990}
L.~A.~Caffarelli,
\newblock \emph{A localization property of viscosity solutions to the Monge--Amp\`ere equation and their strict convexity},
\newblock Ann. of Math. (2) \textbf{131} (1990), no.~1, 129--134.

\bibitem{Caffarelli1990W2p}
L.~A.~Caffarelli,
\newblock \emph{Interior $W^{2,p}$ estimates for solutions of the Monge--Amp\`ere equation},
\newblock Ann. of Math. (2) \textbf{131} (1990), no.~1, 135--150.

\bibitem{Caffarelli1991}
L.~A.~Caffarelli,
\newblock \emph{Some regularity properties of solutions of Monge Amp\`ere equation},
\newblock Comm. Pure Appl. Math. \textbf{44} (1991), no.~8--9, 965--969.

\bibitem{Calabi1958}
E.~Calabi,
\newblock \emph{Improper affine hyperspheres of convex type and a generalization of a theorem by K.~J\"orgens},
\newblock Michigan Math. J. \textbf{5} (1958), 105--126.

\bibitem{CCK2024}
B.~Choi, K.~Choi, and S.~Kim,
\newblock \emph{Translating surfaces under flows by sub-affine-critical powers of Gauss curvature},
\newblock arXiv:\allowbreak 2104.13186v2 (2024).

\bibitem{CCKExistence2025}
B.~Choi, K.~Choi, and S.~Kim,
\newblock \emph{Continuous family of surfaces translating by powers of Gauss curvature},
\newblock Adv. Math. \textbf{483} (2025), Paper No.~110681.

\bibitem{CCKGrowth2025}
B.~Choi, K.~Choi, and S.~Kim,
\newblock \emph{Monge--Amp\`ere equations with right-hand sides of polynomial growth},
\newblock Nonlinear Anal. \textbf{260} (2025), Paper No.~113849.

\bibitem{CW2006}
K.-S.~Chou and X.-J.~Wang,
\newblock \emph{The $L_p$-Minkowski problem and the Minkowski problem in centroaffine geometry},
\newblock Adv. Math. \textbf{205} (2006), 33--83.

\bibitem{CollinsTong2025}
T.~C.~Collins and F.~Tong,
\newblock \emph{Boundary regularity of optimal transport maps on convex domains},
\newblock Acta Math., to appear; arXiv:\allowbreak 2507.05395 (2025).

\bibitem{DaskalopoulosLee2012}
P.~Daskalopoulos and K.-A.~Lee,
\newblock \emph{Fully degenerate Monge--Amp\`ere equations},
\newblock J. Differential Equations \textbf{253} (2012), no.~6, 1665--1691.

\bibitem{daskalopoulos2009monge}
P.~Daskalopoulos and O.~Savin,
\newblock \emph{On Monge--Amp\`ere equations with homogeneous right-hand sides},
\newblock Comm. Pure Appl. Math. \textbf{62} (2009), no.~5, 639--676.

\bibitem{Diaz2015}
G.~D\'iaz and J.~I.~D\'iaz,
\newblock \emph{On the free boundary associated with the stationary Monge--Amp\`ere operator on the set of non strictly convex functions},
\newblock Discrete Contin. Dyn. Syst. \textbf{35} (2015), no.~4, 1447--1468.

\bibitem{DiBenedettoFriedman1986}
E.~DiBenedetto and A.~Friedman,
\newblock \emph{Bubble growth in porous media},
\newblock Indiana Univ. Math. J. \textbf{35} (1986), no.~3, 573--606.

\bibitem{dipierro2023classification}
S.~Dipierro, A.~Karakhanyan, and E.~Valdinoci,
\newblock \emph{Classification of global solutions of a free boundary problem in the plane},
\newblock Interfaces Free Bound. \textbf{25} (2023), no.~3, 455--490.

\bibitem{Dive1931}
P.~Dive,
\newblock \emph{Attraction des ellipso{\"\i}des homog{\`e}nes et r{\'e}ciproques d'un th{\'e}or{\`e}me de Newton},
\newblock Bull. Soc. Math. France \textbf{59} (1931), 128--140.

\bibitem{EberleFigalliWeiss2025}
S.~Eberle, A.~Figalli, and G.~S.~Weiss,
\newblock \emph{Complete classification of global solutions to the obstacle problem},
\newblock Ann. of Math. (2) \textbf{201} (2025), no.~1, 167--224.

\bibitem{EberleShahgholianWeiss2023}
S.~Eberle, H.~Shahgholian, and G.~S.~Weiss,
\newblock \emph{On global solutions of the obstacle problem},
\newblock Duke Math. J. \textbf{172} (2023), no.~11, 2149--2193.

\bibitem{FriedmanSakai1986}
A.~Friedman and M.~Sakai,
\newblock \emph{A characterization of null quadrature domains in $\mathbb R^N$},
\newblock Indiana Univ. Math. J. \textbf{35} (1986), no.~3, 607--610.

\bibitem{GalvezJimenezMira2015}
J.~A.~G\'alvez, A.~Jim\'enez, and P.~Mira,
\newblock \emph{A classification of isolated singularities of elliptic Monge--Amp\`ere equations in dimension two},
\newblock Comm. Pure Appl. Math. \textbf{68} (2015), no.~12, 2085--2107.

\bibitem{GiannopoulosMilman2000}
A.~A.~Giannopoulos and V.~D.~Milman,
\newblock \emph{Extremal problems and isotropic positions of convex bodies},
\newblock Israel J. Math. \textbf{117} (2000), 29--60.

\bibitem{HanHuang2026}
H.~Han and G.~Huang,
\newblock \emph{Regularity of free boundary for Monge--Amp\`ere obstacle problem up to the boundary},
\newblock arXiv:\allowbreak 2610.08351 (2026).

\bibitem{Hug1999}
D.~Hug,
\newblock \emph{Absolute continuity for curvature measures of convex sets II},
\newblock Math. Z. \textbf{232} (1999), 437--485.

\bibitem{HuangTangWang2024}
G.~Huang, L.~Tang, and X.-J.~Wang,
\newblock \emph{Regularity of free boundary for the Monge--Amp\`ere obstacle problem},
\newblock Duke Math. J. \textbf{173} (2024), no.~12, 2259--2313.

\bibitem{JiangTrudingerXu2026}
F.~Jiang, N.~S.~Trudinger, and Q.-Q.~Xu,
\newblock \emph{On a class of Cauchy problems with applications in nonlinear partial differential equations},
\newblock J. Differential Equations \textbf{475} (2026), Paper No.~114525.

\bibitem{JTX2025}
T.~Jin, X.~Tu, and J.~Xiong,
\newblock \emph{Regularity and classification of the free boundary for a Monge--Amp\`ere obstacle problem},
\newblock Adv. Math. \textbf{482} (2025), part B, Paper No.~110602.

\bibitem{JTXSingular2025}
T.~Jin, X.~Tu, and J.~Xiong,
\newblock \emph{On the singular set of the free boundary for a Monge--Amp\`ere obstacle problem},
\newblock Comm. Anal. Geom. \textbf{34} (2026), no.~3, 825--844.

\bibitem{JinXiong2016}
T.~Jin and J.~Xiong,
\newblock \emph{Solutions of some Monge--Amp\`ere equations with isolated and line singularities},
\newblock Adv. Math. \textbf{289} (2016), 114--141.

\bibitem{Jorgens1954}
K.~J\"orgens,
\newblock \emph{\"Uber die L\"osungen der Differentialgleichung $rt-s^2=1$},
\newblock Math. Ann. \textbf{127} (1954), 130--134.

\bibitem{Jorgens1955}
K.~J\"orgens,
\newblock \emph{Harmonische Abbildungen und die Differentialgleichung $rt-s^2=1$},
\newblock Math. Ann. \textbf{129} (1955), 330--344.

\bibitem{Lewy1979}
H.~Lewy,
\newblock \emph{An inversion of the obstacle problem and its explicit solution},
\newblock Ann. Scuola Norm. Sup. Pisa Cl. Sci. (4) \textbf{6} (1979), no.~4, 561--571.

\bibitem{LonginettiSalani2007}
M.~Longinetti and P.~Salani,
\newblock \emph{On the Hessian matrix and Minkowski addition of quasiconvex functions},
\newblock J. Math. Pures Appl. (9) \textbf{88} (2007), no.~3, 276--292.

\bibitem{Lutwak1991}
E.~Lutwak,
\newblock \emph{Extended affine surface area},
\newblock Adv. Math. \textbf{85} (1991), no.~1, 39--68.

\bibitem{Lutwak1993}
E.~Lutwak,
\newblock \emph{The Brunn--Minkowski--Firey theory. I. Mixed volumes and the Minkowski problem},
\newblock J. Differential Geom. \textbf{38} (1993), 131--150.

\bibitem{LutwakYangZhang2004}
E.~Lutwak, D.~Yang, and G.~Zhang,
\newblock \emph{On the $L_p$-Minkowski problem},
\newblock Trans. Amer. Math. Soc. \textbf{356} (2004), no.~11, 4359--4370.

\bibitem{MZ1998}
F.~Merle and H.~Zaag,
\newblock \emph{Optimal estimates for blowup rate and behavior for nonlinear heat equations},
\newblock Comm. Pure Appl. Math. \textbf{51} (1998), no.~2, 139--196.

\bibitem{Milman2025}
E.~Milman,
\newblock \emph{Centro-affine differential geometry and the log-Minkowski problem},
\newblock J. Eur. Math. Soc. \textbf{27} (2025), no.~2, 709--772.

\bibitem{MooneySen2026}
C.~Mooney and A.~Sen,
\newblock \emph{Regularity theory for the equation $\det(D^2u)=f(u)$ in two dimensions},
\newblock arXiv:\allowbreak 2611.06978 (2026).

\bibitem{Pogorelov1972}
A.~V.~Pogorelov,
\newblock \emph{On the improper convex affine hyperspheres},
\newblock Geom. Dedicata \textbf{1} (1972), no.~1, 33--46.

\bibitem{Sakai1981}
M.~Sakai,
\newblock \emph{Null quadrature domains},
\newblock J. Analyse Math. \textbf{40} (1981), 144--154.

\bibitem{Savin2005}
O.~Savin,
\newblock \emph{The obstacle problem for Monge Ampere equation},
\newblock Calc. Var. Partial Differential Equations \textbf{22} (2005), no.~3, 303--320.

\bibitem{SavinYuCones2025}
O.~Savin and H.~Yu,
\newblock \emph{Stable and minimizing cones in the Alt--Phillips problem},
\newblock arXiv:\allowbreak 2502.18192 (2025).

\bibitem{SavinYuConcentration2026}
O.~Savin and H.~Yu,
\newblock \emph{Concentration of cones in the Alt--Phillips problem},
\newblock Arch. Ration. Mech. Anal. \textbf{250} (2026), Paper No.~35.

\bibitem{Schneider2014}
R.~Schneider,
\newblock \emph{Convex Bodies: The Brunn--Minkowski Theory},
\newblock second expanded edition, Encyclopedia of Mathematics and its Applications \textbf{151}, Cambridge University Press, Cambridge, 2014.

\bibitem{wu2022fully}
Y.~Wu and H.~Yu,
\newblock \emph{On the fully nonlinear Alt--Phillips equation},
\newblock Int. Math. Res. Not. IMRN (2022), no.~11, 8540--8570.

\end{thebibliography}
\end{document}